\documentclass[12pt,a4paper]{article}
\usepackage{amsmath,amscd,amssymb,amsthm,amstext}

\theoremstyle{plain}
   \newtheorem{theorem}{Theorem}[section]
   \newtheorem{proposition}[theorem]{Proposition}
   \newtheorem{lemma}[theorem]{Lemma}
   
   \theoremstyle{definition}
   
   \newtheorem{definition}[theorem]{Definition}
   
   \theoremstyle{remark}
   
   \newtheorem{remark}[theorem]{Remark}

\newcommand{\FF}{{\mathbb F}}
\newcommand{\ZZ}{{\mathbb Z}}

\newcommand{\U}{{\mathcal U}}
\newcommand{\K}{{\mathcal K}}

\newcommand{\Alg}{{{\mathcal A}lg}}
\newcommand{\R}{{\mathcal R}}
\newcommand{\sa}{{\mathcal A}}

\newcommand{\Tor}{{\operatorname{Tor}}}

\newcommand{\simp}{_\bullet}

\newcommand{\drf}{{\overline \Omega}}

\newcommand{\trunc}[2]{{T_{#1}(#2)}}
\newcommand{\grtrunc}[3]{{T_{#1 , #2}(#3)}}
\newcommand{\ddr}{{\bf d}}
\newcommand{\gshuf}[2]{{V(#1,#2)}}
\newcommand{\Sq}{\operatorname{Sq}}
\newcommand{\rp}{{\mathbb R} \operatorname{P}}
\newcommand{\cp}{{\mathbb C} \operatorname{P}}
\newcommand{\hp}{{\mathbb H} \operatorname{P}}
\newcommand{\cayley}{{\bf Ca} \operatorname{P}^2} 
\newcommand{\kp}{{\mathbb K} \operatorname{P}}

\newcommand{\thom}[1]{{Th(#1)}}
\newcommand{\rhz}[2]{{\widetilde H_{#1}(#2;\ZZ )}}
\newcommand{\hz}[2]{{H_{#1}(#2;\ZZ )}}
\newcommand{\rcoh}[2]{{\widetilde H^{#1}(#2)}}
\newcommand{\coh}[2]{{H^{#1}(#2)}}

\newcommand{\num}[1]{{\mathbb #1}}
\newcommand{\numcayley}{{\bf Ca}}
\newcommand{\lm}{{CT}}

\title{A splitting result for the free loop space of 
spheres and projective spaces}
\author{Marcel B\" okstedt \& Iver Ottosen}
\date{May 17, 2004}

\begin{document}
\maketitle

\begin{abstract}
Let $X$ be a 1-connected compact space such that the algebra 
$H^*(X;\FF_2 )$ is generated by one single element. 
We compute the cohomology of the free loop space 
$H^*(\Lambda X;\FF_2 )$ including the Steenrod algebra action. 
When $X$ is a projective space $\cp^n$, $\hp^n$, the Cayley 
projective plane $\cayley$ or a sphere $S^m$ we obtain a 
splitting result for integral and mod two cohomology of the 
suspension spectrum $\Sigma^\infty (\Lambda X)_+$. The splitting
is in terms of $\Sigma^\infty X_+$ and the Thom spaces 
$\thom {q\tau}$, $q\geq 0$ of the $q$-fold Whitney sums of 
the tangent bundle $\tau$ over $X$.
\end{abstract}

\section{Introduction}

The main part of this paper is a computation of 
mod two cohomology of the free loop space on a space $X$,
with the property that 
$H^*(X,\FF_2)$ is a truncated polynomial algebra. 

This calculation has been done by various other methods before,
see \cite{Z}, \cite{KY} and \cite{CC}. We want to 
make four points in order to explain what is new 
in this paper.  

First, we are using a new general method.
The idea is to resolve the space $X$ as a degreewise free
cosimplicial space, and then apply the functor ``free
loop space'' degreewise. The Bousfield spectral sequence
associated to this cosimplicial space gives us a spectral 
sequence converging to the (co)-homology
of the free loop space of $X$, at least if $X$ is
simply connected. 

This spectral sequence is a generalization of the
Eilenberg-Moore spectral sequence.
Eilenberg-Moore spectral sequences has been applied in this situation
before.
The main difference from previous work is that the cosimplicial space 
we consider is not the same as
the cosimplicial spaces obtained from homotopy pull back squares like
\[
\begin{CD}
\Lambda X @>>> X \\
@VVV      @V{\Delta}VV \\
X @>{\Delta}>> X\times X.  
\end{CD}
\]
The cosimplicial space we use is related to the cosimplicial
resolution of $X$ (in the sense of \cite{BK}).
In this paper we are using
this method to treat some serious cases. From this point of view, we
are dealing with test cases for the spectral sequence.

Secondly, the spectral sequence does give us information about the
action of the Steenrod algebra on the cohomology of the free
loop space. We determine this action for
spaces belonging to the following list : 
$X=S^n$, $X=\cp^n$, $X=\hp^n$ and $X=\cayley$.
Except for the easy sphere case, these results are new.

Thirdly, the circle group acts in an obvious way on a free
loop space, and one can consider the Borel construction of this
action. 
The method we consider extends in a natural way 
(\cite{BO2}) to a spectral sequence converging to
the Borel construction. The structures discovered in
\cite{CS} has recently boosted the interest in the homology 
of the Borel construction.
In spite of this, this homology has not
been calculated. We intend to study the corresponding spectral 
sequence. From this point of view, this paper is
a preliminary result.

Fourthly, for the spaces in the short list above, 
we obtain a surprising result.
The spaces are all manifolds, and even symmetric spaces. In each case
we consider the tangent bundle $\tau$. 
Let $q\tau$ be the $q$-fold Whitney sum of the tangent bundle with
itself. We can consider $q\tau$ as a sub bundle of $(q+1)\tau$.
This induces an inclusion of the corresponding Thom spaces,
and we define the space
$C_q(X)$ as the cofiber of the cofibration
\[
\thom {q\tau} \to \thom {(q+1)\tau}.
\]

The result is that cohomology of $\Sigma^\infty (\Lambda X_+)$ 
with coefficients in $\FF_2$ is the same $\sa$-module 
as the cohomology of
\[
\Sigma^\infty \Big(X_+ \vee \bigvee_{q\geq 0} \Sigma^{(m-2)(q+1)}
C_q(X)\Big),
\]
where $m$ denotes the degree of the algebra generator for 
$H^*(X;\FF_2)$.

A splitting discovered by Ziller \cite{Z} shows that the homology 
with $\ZZ$-coefficients of the two spectra also agree as abelian 
groups.

This strongly suggest the possibility that this wedge is
actually homotopy equivalent to $\Sigma^\infty (\Lambda X)_+$.
There is some additional evidence for such a splitting. 
If $X=S^n$ for $n\geq 2$, the spectrum $\Sigma^\infty (\Lambda X_+)$
is known and simple. It is a wedge of sphere spectra and 
Moore spectra for the group $\ZZ /2$. 
In fact, in this special case, we do obtain a stable splitting of
spectra as above.

\begin{remark}
\label{rem1}
There are a few simply connected spaces besides the spheres and 
the projective spaces considered above that have mod two cohomology 
rings that are truncated polynomial algebras. But there are certain 
conditions which the cohomology ring of such a space has to satisfy.
Let $X$ be a 1-connected and compact space. 
Assume that $H^*(X;\FF_2)= \FF_2 [x]/(x^{n+1})$ for some $n\geq 1$  
and let $m$ be the degree of $x$.
If $n\geq 2$ it is known known \cite{A}, \cite{T} that
$m=2,4$ or $m=8$ and $n=2$.
For $n=1$ we have that $m$ can be any integer bigger that $2$.
\end{remark}

{\bf Notation:} All cohomology groups are with $\FF_2$-coefficients
unless stated otherwise. We write $\sa$ for the mod two Steenrod
algebra, $\U$ for the category of unstable $\sa$-modules and
$\K$ for the category of unstable $\sa$-algebras. The degree
of an element $x$ is denoted $|x|$. The category
of non negatively graded $\FF_2$-algebras with the property
that $a^2=a$ if $|a|=0$ is denoted $\Alg$.
Note that we have a forgetful functor $\K \to \Alg$. 

We write $\Omega : \Alg \to \Alg$ for the de Rham complex functor.
Thus $\Omega (A)$ for an $A\in \Alg$ means the
free $A$-algebra on generators $\{ \ddr x| x\in A \}$ of degree 
$|\ddr x|=|x|-1$ modulo the relations 
$\ddr (x+y)= \ddr x + \ddr y$,
$\ddr (xy) = \ddr (x)y+x\ddr (y)$ and $(\ddr x)^2=0$. 
We always write $\ddr$ for the
de Rham differential in order to distinguish it from the 
simplicial face maps $d_i$.
By $s\mathcal C$ we mean the category of simplicial objects 
in the category $\mathcal C$.

\section{Results via a simplicial resolution in $s\Alg$}

In this section we compute the $E_2$-term of a spectral
sequence form \cite{BO2} converging toward $H^*(\Lambda X)$ in the
case where $X$ is a 1-connected, compact space and $H^*X$ a truncated 
polynomial algebra $\FF_2 [x]/(x^{n+1})$. 
The $E_2$-term is given by the (non abelian) derived functors of a 
twisted de Rham complex functor $\drf$ as follows
$$E_2^{-m,t} = H_m(H^*X ; \drf)^t := (L_m\drf)(H^*X)^t.$$ 
See \cite{BO2} Definition 7.1 and sections 5, 6 for the precise
definitions.

Remark \ref{rem1} shows that the operation 
$\lambda (x)=\Sq^{|x|-1}x$ vanish on $H^*X$. 
By \cite{BO2} Theorem 6.4, Theorem 7.6 and Proposition 7.3
this implies that 
$$H_*(H^*X;\drf ) \cong H_*(H^*X;\Omega ),$$
so it suffices to compute the derived functors of the ordinary
de Rham complex functor over the simplicial category $s\Alg$. 
It turns out that the two cases $n$ even, $n$ odd are essentially 
different. The result of the computation appears in Theorem \ref{main1}.

Let $\trunc n x$ denotes the graded commutative $\FF_2$-algebra 
$\FF_2[x]/(x^{n+1})$ where $n\geq 1$ and $|x|\geq 1$.

\begin{theorem}
There is an almost free simplicial resolution $R\simp \in s\Alg$ of 
$\trunc n x$ as follows: 
$R_q=\FF_2 [x,y_1,\dots, y_q]$ for $q\geq 0$ where $|y_i|=(n+1)|x|$ 
and $s_i(x)=x$, $d_i(x)=x$ for all $i$ and 
\begin{align*}
& s_i(y_j) = 
\begin{cases}
y_j &, i \geq j \\
y_{j+1} &, i<j 
\end{cases} \\
& d_i(y_j) = 
\begin{cases} 
x^{n+1} &, i=0, j=1 \\
y_{j-1} &, i<j, j>1 \\
y_j &, i\geq j , j<q \\
0 &, i=q, j=q.
\end{cases}
\end{align*}
\end{theorem}

\begin{proof}
Similar to the proof of Proposition 8.2 in \cite{BO2}.
\end{proof}

For a simplicial $\FF_2$-vector space $V$ we write
$C_*(V)$ for the associated chain complex with $C_q(V)=V_q$ and 
differential $\sum_{i=0}^q d_i$. We write $N_*(V)$ for the normalized 
chain complex with $N_q(V)=\cap_{1\leq i \leq q} \ker (d_i)$ and 
differential $d_0$. Finally, we write $D_*(V)$ for the sub complex 
of $C_*(V)$ generated by the degenerate elements \cite{ML}.

\begin{definition}
The element $\omega_q \in \Omega (R_q)$ is defined by
$\omega_q = \ddr y_1 \dots \ddr y_q$ for $q>0$ and $\omega_0=1$.
We use the notation
\begin{align*}
& (\omega_q)_i = \ddr y_1 \dots \widehat{ \ddr y_i} \dots \ddr y_q
& \text{ for } & 1\leq i \leq q, \\
& (\omega_q)_{i,j} = \ddr y_1 \dots \widehat{ \ddr y_i} \dots 
\widehat{ \ddr y_j} \dots \ddr y_q
& \text{ for } & 1\leq i <j \leq q.
\end{align*}
Define $\alpha_q ,\beta_q \in \Omega (R_q)$ by
$$
\alpha_q = \ddr x \omega_q \quad , \quad
\beta_q = x\omega_q + \ddr x \sum_{i=1}^q y_i(\omega_q)_i \quad ,
$$ 
for $q>0$ and $\alpha_0=\ddr x$, $\beta_0 = x$.
\end{definition}

\begin{lemma}
\label{cycle_deg}
For $q\geq 0$ one has that  
$\omega_q$ is a cycle in $N_q\Omega (R\simp )$ for $n$ odd and
$\alpha_q$ , $\beta_q$ are cycles in $ N_q \Omega (R\simp )$ 
for $n$ even. Furthermore, the elements
$(\omega_{q+1} )_r$, $(\omega_{q+1})_{j,k}$, 
$y_r (\omega_{q+1} )_{j,k}$ lie in $D_{q+1}\Omega (R\simp )$
for any $n$ where $1\leq j<k \leq q+1$ and $1\leq r \leq q+1$.
\end{lemma}

\begin{proof}
For $n$ odd we have $d_0\omega_q=0$ since 
$d_0 \ddr y_1 = \ddr (x^{n+1})=0$,
$d_i \omega_q =0$ for $0<i<q$ since $(\ddr y_i)^2=0$ and 
$d_q\omega_q =0$ since $d_q \ddr y_q =0$.
For $n$ even we have $d_i\beta_q =0$ for all $i$ since
$$
d_i(x\omega_q )= 
\begin{cases}
x^{n+1}\ddr x\omega_{q-1} &, i=0 \\
0 &, i>0,
\end{cases}
$$
$$
d_i(y_j\ddr x(\omega_q )_j) = 
\begin{cases}
x^{n+1}\ddr x\omega_{q-1} &, i=0;j=1 \\
y_i\ddr x\omega_{q-1} &, 0<i<q; j=i,i+1 \\
0 &, \text{otherwise.}
\end{cases}
$$
Note that $\ddr x\beta_q =x\alpha_q$ such that $xd_i\alpha_q =0$ and
hence $d_i\alpha_q =0$.

The last statement follows directly by the formulas for the degeneracy maps.
\end{proof}

By a $(p,q)$-shuffle $(\mu , \nu)$ 
we understand two ordered subsequences
$\mu_1 < \dots <\mu_p$ and $\nu_1 <\dots <\nu_q$ 
of the set
$\{ 0 , 1, \dots ,p+q-1\}$
such that
$$\{ \mu_1 , \dots ,\mu_p\} \sqcup \{ \nu_1 ,\dots ,\nu_q \} =
\{ 0, 1, \dots ,p+q-1 \} .$$
Furthermore, we use the notation $s_\mu = s_{\mu_p} \dots s_{\mu_1}$.

\begin{lemma}
\label{shuffle}
Assume $(\mu, \nu )$ is a $(p,q)$-shuffle and $j\geq 1$.
Then $s_\mu (y_j)=y_{\nu_j+1}$ and hence
$$s_\mu (\omega_q )= \prod_{t\in \{ \nu_1 ,\dots ,\nu_q \} }
\ddr y_{t+1} \quad , \quad
s_\mu ((\omega_q )_j)=
\prod_{t\in \{ \nu_1 ,\dots ,\widehat{\nu_j}, \dots ,\nu_q \} }
\ddr y_{t+1}.
$$
\end{lemma}

\begin{proof}
We compute $s_\mu (y_j)$ by the formulas for the degeneracy maps.
If $\mu_1 \geq j$ then $\nu_j = j-1$ and $s_\mu (y_j)=y_j$.
If $\mu_2 \geq j+1$ and $\mu_1 <j$ then $\nu_j=j$ and 
$s_\mu(y_j)=y_{j+1}$. In general, if $\mu_{r+1} \geq j+r$ and
$\mu_r < j+r-1$ then $\nu_j = j+r-1$ and $s_\mu (y_j)=y_{j+r}$.
\end{proof}

By results of Bousfield and Dwyer \cite{D} we have that
$H_*(\trunc n x;\Omega ) = \pi_* \Omega (R\simp )$ is
a graded $\FF_2$-algebra equipped with homotopy operations
$$\delta_i: H_q(\trunc n x; \Omega ) \to H_{q+i}(\trunc n x ; \Omega ) 
\quad , \quad 2\leq i \leq q.$$ 
The following result gives a complete determination of this structure.

\begin{theorem}
\label{main1}
Let $n$ be a positive integer.
For $n$ odd, there is an isomorphism of algebras
$$H_*(\trunc n x ; \Omega ) \cong 
\trunc n x \otimes \Lambda (\ddr x) \otimes \Gamma [\omega ]$$
where $x, \ddr x \in H_0$ and $\gamma_q(\omega )\in H_q$ 
The homotopy operations are given by 
$$\delta_i (\gamma_q(\omega)) = 
\begin{cases}
\binom {2q-1} q \gamma_{2q}(\omega ) &, i=q \\
0 &, i<q.
\end{cases}
$$
The algebra generators are represented 
by cycles in the normalized chain complex $(N_*\Omega R,d_0)$ with
$N_q\Omega R= \cap_{i>0} \ker d_i$ as follows:
$x=[x]$, $\ddr x=[\ddr x]$ and $\gamma_i (\omega ) =[\omega_i]$.

For $n$ even, we have 
$$H_*(\trunc n x;\Omega ) \cong 
\FF_2 \oplus \FF_2 [x]/(x^n) <a_q,b_q|q\geq 0>$$
where $x\in H_0$ and $a_q,b_q\in H_q$.
The algebra structure is given by the relations
$$a_pa_q=0,\quad b_pb_q = \binom {p+q} p xb_{p+q}, \quad
a_pb_q = \binom {p+q} p xa_{p+q}.$$ 
The homotopy operations are given by $\delta_i (a_q)=0$ and
$$\delta_i (b_q) = 
\begin{cases}
\binom {2q-1} q xb_{2q} &, i=q \\
0 &, i<q.
\end{cases}
$$
The generators are represented by cycles as follows: 
$x=[x]$, $a_q=[\alpha_q]$, $b_q=[\beta_q]$.
\end{theorem}

\begin{remark} 
The description in the theorem implies that 
$H_0(\trunc n x;\Omega )\cong \Omega (\trunc n x)$.
Note also that
$$\binom {2q-1} q = \begin{cases} 1 &, q =2^r,r \geq 0 \\ 0 &, 
\text{ otherwise.} \end{cases}$$
\end{remark}

\begin{proof}
We start by computing the derived functors additively.
The algebra $\trunc n x$ is the pushout of the diagram
$\FF_2 \leftarrow \FF_2 [y] \rightarrow \FF_2 [x]$ where 
$y\mapsto x^{n+1}$. Note that $\FF_2 [x]$ is free as a module over
$\FF_2 [y]$. By \cite{BO2} Proposition 6.3 there is a Quillen spectral
sequence as follows:
$$E_{i,j}^2=\Tor_i^{H_*(\FF_2 [y];\Omega )}
(\FF_2 , H_*(\FF_2[x];\Omega ))_j 
\Rightarrow H_{i+j}(\trunc n x;\Omega ).$$
The algebras $\FF_2 [x]$ and $\FF_2 [y]$ are free objects in $\Alg$ so 
we find that $E_{i,j}^2=0$ for $j>0$ and
$$H_i(\trunc n x;\Omega )\cong E_{i,0}^2 \cong 
\Tor_i^{\Omega (\FF_2 [y])}(\FF_2 , \Omega (\FF_2 [x])).$$

There is a resolution $(K_*, \partial )$ of $\FF_2$ by free
$\Omega (\FF_2 [y])$-modules as follows:
$$K_* = \Lambda (v) \otimes \Gamma [w] \otimes \Omega (\FF_2 [y]), 
\quad \partial v = y, 
\quad \partial \gamma_i (w)=\gamma_{i-1} (w)dy$$
where $v\in K_1$ and $\gamma_i (w) \in K_i$. We tensor $K_*$ 
with $\Omega (\FF_2 [x])$ over $\Omega (\FF_2 [y])$ and
get a complex $(C_*, \partial )$ with 
$$C_*=\Lambda (V)\otimes \Gamma [w]\otimes \Omega (\FF_2 [x]), \quad 
\partial v = x^{n+1}, \quad 
\partial \gamma_i (w)=(n+1)\gamma_{i-1}(w)x^ndx$$
By computing the homology of $(C_*,\partial )$ we find that for $n$
even $H_0(\trunc n x;\Omega )= \Omega (\trunc n x)$ and 
$$H_q(\trunc n x; \Omega )= 
\FF_2 [x]/(x^n) <dx \gamma_q (w),x\gamma_q (w)+dxv\gamma_{q-1}(w)>,
\quad q>0.$$
For $n$ odd we find
$$H_*(\trunc n x;\Omega ) \cong \Omega (\trunc n x)\otimes \Gamma [w].$$

We now show that there is an additive basis for 
$H_q(\trunc n x; \Omega )$ which is represented by cycles in 
the normalized complex as stated. 
Lemma \ref{cycle_deg} shows that these representatives are indeed cycles.
It suffices to show that the associated classes are linearly 
independent since, in each degree, the number of classes agree with 
the dimension in that degree. 

We introduce two gradings on $\Omega R\simp$ preserved by $d_i$.
The wedge grading $w$ which equals the number of wedge factors and
the polynomial grading $p$ with $p(x)=2$, $p(dx)=1$, $p(y_i)=2n+2$,
$p(dy_i)=2n+1$ extended multiplicatively.
Let $\Omega^{r,s}(R_q)$ denote the subspace of elements
$\omega \in \Omega (R_q)$ with $w(\omega )=r$ and $p(\omega )=s$.
We have defined a splitting of simplicial $\FF_2$-vector spaces so it
follows that
$$H_q(\trunc n x;\Omega ) \cong \bigoplus_{r,s\geq 0} 
H_q(\trunc n x;\Omega^{r,s} ).$$

Assume that $n$ is odd. We must show that for fixed $q$ the elements
$$x^t(\ddr x)^\epsilon \omega_q, 
\quad 0\leq t \leq n, \epsilon \in \{ 0 ,1\} $$
are linearly independent. The bidegree is given by 
$(q+\epsilon, 2t+\epsilon +q(2n+1))$. The elements
have different bigradings and it suffices to check that
they individually do not represent zero.
Recall that the composite
$N_*(V) \to C_*(V) \to C_*(V) /D_*(V)$
is an isomorphism for any simplicial $\FF_2$-vector space $V\simp$ 
so we have that
$N_*\Omega (R\simp ) \cap D_*\Omega (R\simp )=0$.

We first check that $[x^t\ddr x \omega_q]\neq 0$. The chains
which might hit the representative lie in 
$\Omega^{q+1,2t+1+q(2n+1)}(R_{q+1})$ which is spanned by the elements
$x^t\ddr x (\omega_{q+1})_j$ for $1\leq j \leq q+1$ together with 
$\omega_{q+1}$ if $t=n$.
By Lemma \ref{cycle_deg} all these elements are degenerate except for 
$\omega_{q+1}$. It follows that the only non trivial normalized chain
in this bigrading is $\omega_{q+1}$. But $d_0(\ddr y_1)=0$ such that 
$d_0\omega_{q+1} = 0$. 

We then check that $[x^t\omega_q ]\neq 0$. The chains which might
hit the representative lie in $\Omega^{q,2t+q(2n+1)}(R_{q+1})$
which is spanned by the elements
$x^t (\omega_{q+1} )_r$, $x^{t+n}\ddr x (\omega_{q+1} )_{j,k}$
together with  $x^{t-1}\ddr x y_r(\omega_{q+1} )_{j,k}$ if $t>0$
where $1\leq r \leq q+1$, $1\leq j<k\leq q+1$.
By Lemma \ref{cycle_deg} all these elements are degenerate.
So there are no normalized chains that can have $x^t \omega_q$
as boundary.

Assume that $n$ is even. We must show that the elements
$$x^t\alpha_q , x^t\beta_q, \quad 1\leq t \leq n-1$$
are linearly independent. The bidegrees are
$$||x^t\alpha_q ||=(q+1,2t+1+q(2n+1)) \quad , \quad 
||x^t\beta_q ||=(q,2t+2+q(2n+1))$$
so the elements have different bigradings and it suffices to check that
they individually do not represent zero.

We first check that $[x^t\alpha_q ]\neq 0$. The chains
which might hit the representative lie in
$\Omega^{q+1,2t+1+q(2n+1)}(R_{q+1})$ which is spanned by the elements
$x^{t+1}\ddr x (\omega_{q+1} )_r$ for $1\leq r \leq q+1$ 
By Lemma \ref{cycle_deg} all these elements are degenerate.

We then check that $[x^t\beta_q ]\neq 0$. The chains which might
hit the representative lie in $\Omega^{q,2t+2+q(2n+1)}(R_{q+1})$
which is spanned by the elements
$$
x^{t+1} (\omega_{q+1} )_r, \quad
x^{t+1+n}\ddr x (\omega_{q+1} )_{j,k}, \quad
x^t\ddr x y_r(\omega_{q+1} )_{j,k}, 
$$
for $1\leq r\leq q+1$, $ 1\leq j<k\leq q+1$. 
By Lemma \ref{cycle_deg} all these elements are degenerate.

We now prove that the algebra structure and homotopy operations
are as stated. The product comes from the chain map
$$\begin{CD}
\rho : C_*(\Omega R\simp ) \otimes C_*(\Omega R\simp ) @>g>>
C_*(\Omega R\simp \otimes \Omega R\simp ) @>C_*(\mu )>> 
C_*(\Omega R\simp ),
\end{CD}$$
where $g$ denotes the Eilenberg-MacLane map \cite{ML} VIII.8. 
We use the results from \cite{BO3} to find the homotopy operations
$$\delta_i : H_q(\trunc n x ;\Omega ) \to H_{q+i}(\trunc n x ;\Omega ), \quad 2\leq
i \leq q.$$
For a cycle $z\in N_q\Omega (R\simp )$ we have 
$\delta_i ([z])=[\delta_i (z)]$ where 
$$\delta_i (z) = \sum_{(\mu ,\nu ) \in \gshuf i q} s_{\nu}(z)s_{\mu}(z).$$

Assume that $n$ is odd. By the formulas for the shuffle map 
and Lemma \ref{shuffle} we find that
$$
\rho (\omega_p \otimes \omega_q ) = \sum_{(\mu ,\nu )} 
s_{\nu} (\omega_p )
s_{\mu} (\omega_q )
= \sum_{(\mu ,\nu )} \omega_{p+q} = \binom {p+q} p \omega_{p+q} ,
$$
where the sums are over all $(p,q)$ shuffles.
Thus the algebra structure is as stated.

The top homotopy operation is as stated since by 
Lemma \ref{shuffle} we have

$$
\delta_q (\omega_q ) 
= \sum_{(\mu ,\nu ) \in \gshuf q q} s_\nu (\omega_q )s_\mu (\omega_q )
= \sum_{(\mu ,\nu ), \mu_1=0} \omega_{2p} 
= \binom {2p-1} {p-1}\omega_{2p}.
$$

For the other homotopy operations we note that when 
$(\mu ,\nu) \in \gshuf i q$ for $2\leq i \leq q-1$
then both $s_\nu (\omega_q )$ and $s_\mu (\omega_q )$ contains
a factor $\ddr y_1$ by Lemma \ref{shuffle} so the operations vanish
since $(\ddr y_1)^2=0$.

Assume that $n$ is even. We have that
$\rho (\alpha_p \otimes \alpha_q)=0$. So the first product
formula holds. By the formulas defining $\beta_p$ and
$\beta_q$ we find that
\begin{align*}
\rho (\beta_p \otimes \beta_q) &= x^2 \rho (\omega_p \otimes \omega_q) \\
&+ x\ddr x \big( 
\sum_{i=1}^p \rho (y_i (\omega_p )_i \otimes \omega_q)
+\sum_{j=1}^q \rho (\omega_p \otimes y_j (\omega_q )_j)
\big).
\end{align*}

The first term equals $\binom {p+q} p x^2 \omega_{p+q}$. 
By Lemma \ref{shuffle} we find that
$$\rho (\omega_p \otimes y_j(\omega_q )_j)=
\sum_{(\mu ,\nu )} s_\nu (\omega_p ) s_\mu (y_j) s_\mu ((\omega_q )_j) =
\sum_{(\mu ,\nu )} y_{\nu_j+1} (\omega_{p+q})_{\nu_j+1} ,$$
and by symmetry
$$\rho (y_i (\omega_p )_i \otimes \omega_q)=
\sum_{(\mu ,\nu )} y_{\mu_i+1} (\omega_{p+q})_{\mu_i+1}.$$
We insert these above and find that
$\rho (\beta_p \otimes \beta_q) = \binom {p+q} p x\beta_{p+q}$.
This gives the second formula for the product. The third
follows by multiplication with $\ddr x$ on the chain level.

A similar argument shows that 
$\delta_q (\beta_q )=\binom {2q-1} {q-1} x\beta_{2q}$. 
For $2\leq i \leq q-1$ we get $\delta_i (\beta_q )=0$ since 
$(\ddr y_1)^2=0$. Finally, we have that $\delta_k (\alpha_q )=0$
for all $k$ since $(\ddr x)^2=0$. 
\end{proof}

\section{Normalized chain complex manipulations}
\label{section:normalized}

This very technical section is needed in the discussion of the
Steenrod algebra action in the next section.

Let $A\simp$ be a simplicial ring. In a DGA it follows 
directly from the definition that the product of two cycles 
is a cycle, and that the product of a cycle with a boundary 
is a boundary. The Eilenberg-Zilber equivalence allows you to
transfer this statement to simplicial rings. 
We will only use a rather special case.

\begin{definition}
Let $m(a,b)=s_0(a)s_1(b)-s_0(b)s_1(a)\in A_{n+1}$ 
and $q(a)=s_0(a)s_1(a)\in A_{n+1}$ for $a,b \in A_n$.
\end{definition}

When $n=1$ the map $m$ can be  described as the shuffle map
followed by the multiplication map 
$A_1\otimes A_1 \overset g \to A_2\otimes A_2 
\overset \mu \to A_2.$

The bilinear map $m$ is skew-symmetric and $m(a,a)=0$.
In case $A\simp$ is a simplicial $\FF_2$-algebra we have that
$q(a+b)=q(a)+q(b)+m(a,b)$ so the two constructions formally behave 
like a symmetric bilinear form and the quadratic form belonging to it. 

We will show that $m$ preserves the normalized chain complex, and
pass to well-defined operations on its homology.
The situation for $q$ is similar, but slightly more complicated.

\begin{lemma}
\label{lemma:products}
Assume that $a,b\in N_n(A)$ and $c\in A_{n-1}$ (in particular 
we may take $c=1$ when $A\simp$ is unital). Then
\begin{itemize}
\item $s_1s_0(c)m(a,b) \in N_{n+1}(A)$.
\item If $d_0(a)c=0=d_0(b)c$, then $s_1s_0(c)m(a,b)$ is a cycle.
\item If $d_0(a)c=0$ is a cycle and $b$ is a boundary, then 
$s_1s_0(c)m(a,b)$ is a boundary.
\end{itemize}
\end{lemma}

\begin{proof}
We compute the face maps on the element $z=s_0s_0(c)m(a,b)$.
\begin{align*}
d_0(z) &= s_0(c)(as_0d_0(b)-s_0d_0(a)b) 
= as_0(cd_0(b))-s_0(cd_0(a))b,\\
d_1(z) &= s_0(c)(ab-ab)=0,\\
d_2(z) &= s_0s_0d_0(c)[s_0d_1(a)b-as_0d_1(b)]=0,\\
d_i(z) &=
s_0s_0d_{i-2}(c)[s_0d_{i-1}(a)s_1d_{i-1}(b)-
s_1d_{i-1}(a)s_0d_{i-1}(b)]=0 \text{ for } i\geq 3.
\end{align*}
The first two statements follows from this.

Now assume that $d_0(a)c=0$ and that $b=d_0(x)$ for some 
$x\in N_{n+1}(A)$. We define
$$y=s_2s_1s_0(c)[s_2s_1(a)s_0(x)-s_2s_0(a)s_1(x)+
s_1s_0(a)s_2(x)]\in A_{n+2}.$$
Note that this element is derived from the shuffle map. 
We compute face maps:
\begin{align*}
d_0(y)&=s_1s_0(c)[s_1s_0d_0(a)x-s_1(a)s_0d_0(x)+s_0(a)s_1d_0(x)]\\
&=s_1s_0(c)s_1s_0d_0(a)x+s_1s_0(c)[s_0(a)s_1(b)-s_1(a)s_0(b)]\\
&=s_1s_0(cd_0(a))x+s_1s_0(c)m(a,b)=s_1s_0(c)m(a,b),\\
d_1(y)&=s_1s_0(c)[s_1(a)x-s_1(a)x+s_0(a)s_1d_1(x)]=0,\\  
d_2(y)&=s_1s_0(c)[s_1(a)s_0d_1(x)-s_0(a)x+s_0(a)x]=0,\\
d_i(y)&=s_1s_0d_{i-2}s_0(c)[s_2s_1d_{i-2}(a)s_0d_{i-1}(x)\\
&-s_2s_0d_{i-2}(a)s_1d_{i-1}(x)+s_1s_0d_{i-2}(a)s_2d_{i-1}(x)]=0
\text{ for } i\geq 3.
\end{align*}
\end{proof}

\begin{lemma}
\label{lemma:squares}
Assume that $a\in N_n(A)$, $b\in N_{n+1}(A)$ and $c\in A_n$ 
(in particular we may take $c=1$ if $A\simp$ is unital).
\begin{itemize}
\item If $ca^2=0$, then $s_0(c)q(a) \in N_{n+1}(A)$.
\item If in addition $cas_0d_0(a)=0$ 
(in particular, if $a$ is a cycle), then $s_0(c)q(a)$ is a cycle.
\item If $s_0(c)b^2=0$, then $s_0(c)q(d_0(b))$ is a boundary.
\end{itemize}  
\end{lemma}

\begin{proof}
Put $y=s_0(c)q(a)=s_0(c)s_0(a)s_1(a)$. 
The first two statements follows by the computation
\begin{align*}
d_0(y)&=cas_0d_0(a),\\
d_1(y)&=ca^2=0,\\
d_2(y)&=s_0d_1(c)s_0d_1(a)a=0,\\
d_i(y)&=s_0d_{i-1}(c)s_0d_{i-1}(a)s_1d_{i-1}(a)=0 
\text{ if } i\geq 3.\\
\end{align*}

We claim that if $s_0(c)b^2=0$, then $s_1s_0(c)s_1(b)s_2(b)$ 
is a chain with boundary $s_0(c)q(d_0(b))$. The claim follows 
by the computation:
$$
d_i(s_1s_0(c)s_1(b)s_2(b))=
\begin{cases}
  s_0(c)s_0d_0(b)s_1d_0(b)=s_0(c)q(d_0(b)) &\text{ if } i=0,\\
  s_0(c)bs_1d_1(b)=s_0(c)b\cdot 0=0 &\text{ if } i=1,\\
  s_0(c)b^2=0 &\text{ if } i=2,\\
  s_1s_0d_{i-2}(c)s_1d_{i-1}(b)s_2d_{i-1}(b)=0 &\text{ if } i\geq 3.
\end{cases}
$$
\end{proof}

\section{The Steenrod algebra action}

The algebra $\trunc n x$ is sometimes an object in $\K$.
For $n\geq 2$ this happens if and only if $|x|$ is a power 
of two. For $n=1$ it happens if and only if $|x|\geq 1$. 
There is at most one $\sa$-action which makes $\trunc n x$
an object in $\K$. 

For any $K \in \K$ there is an $\sa$-action on
$H_*(K; \drf)$. In this section we determine this
action when $K=\trunc n x \in \K$.

We first consider the case where $|x|$ is a power of two.
Define $\grtrunc n t x \in K$ for $n\geq 1$, $t\geq 0$ 
by $\grtrunc n t x = \FF_2[x]/(x^{n+1})$ with $m=|x|=2^t$. It is
an object in $\K$ by the formula
$$\Sq^i (x^j)= \binom {mj} i x^{j+j/m}.$$ 

Let $\Phi: \K \to \K$ denote the doubling functor with
$\Phi (K)^{2j} = K^j$ and $\Phi (K)^{2j+1}=0$ for $K\in \K$.
The operation $\Sq ^{2i}$ act on $\Phi (K)$ via the action of 
$\Sq ^i$ on $K$. The operation $\Sq^{2i+1}$ acts trivially 
on $\Phi (K)$. 
By \cite{BO2} Theorem 6.4 and Theorem 7.6 we have the 
following result:

\begin{proposition}
\label{doubling}
Let $A \in \K$ and assume that $A=\Phi^t (K)$ for some $t\geq 1$ and
$K\in \K$. Then there is an isomorphism of $\sa$-modules for all $q$:
$$H_q(A;\drf ) \cong \pi_q \Omega \Phi^t (B\simp )$$
where $B\simp \in s\K$ is an almost free resolution of $K$.
\end{proposition}

Note that $\grtrunc n t x = \Phi^t \grtrunc n 0 x$ where the $\sa$-action
on $\grtrunc n 0 x$ is given by
\begin{equation}
\label{square}
\Sq^j x^i=\binom i j x^{i+j}.
\end{equation}
Let $R\simp^{\K} \in s\K$ denote a simplicial resolution of 
$\grtrunc n 0 x \in \K$. There is a simplicial algebra map 
$R\simp \to R\simp^{\K}$ which induces a homotopy equivalence. 
We will abuse the notation, and denote the images in $R_q^{\K}$ 
of the classes $x,y_1,\dots,y_q\in R_q$ by the same names.

By the formulas for the boundary maps and (\ref{square}) we 
see that the element 
$\Sq^i(y_1)+\binom {n+1} i x^iy_1\in R_1^{\K}$ 
is a cycle for $1\leq i \leq n+1$.
Since $R\simp^\K$ is a resolution of a discrete ring, 
$\pi_k(R\simp )=0$ for $k>0$ and 
we can find elements $z^{(i)}\in R_2^{\K}$ for $1\leq i \leq n+1$
such that  
$$
d_0(z^{(i)})=\Sq^i(y_1)+\binom {n+1} i x^iy_1,
\quad d_1(z^{(i)})=d_2(z^{(i)})=0.
$$     
As we shall see below this information on $R^{\K}\simp$ turns out
to be sufficient for the computation of the $\sa$-action
on $H_*(\grtrunc n t x;\drf )$.

\begin{theorem}
\label{sqodd}
Let $n\geq 1$, $t\geq 1$ and assume that $n$ is odd. 
Then the $\sa$-action on  
$$H_*(\grtrunc n t x; \drf )\cong 
\grtrunc n t x \otimes \Lambda (\ddr x)\otimes \Gamma [\omega ]$$ 
is given by the following formulas where $i,q\geq 0$:
\begin{align}
\label{sqodd1}
& \Sq^{2^t i} (\gamma_q (\omega ) ) = 
\binom {q(n+1)} i x^i \gamma_q (\omega ) \\
\label{sqodd2}
& \Sq^{2^t i} (x^j) = \binom j i x^{i+j} \\
\label{sqodd3}
& \Sq^{2^t i} (\ddr x) = 
\begin{cases}
\ddr x &, i=0 \\
0 &, i>0.
\end{cases}
\end{align} 
If $k \neq 0$ mod $2^t$ then $\Sq ^k(z)=0$ for all 
$z\in H_*(\grtrunc n t x ;\drf )$.
\end{theorem}

\begin{proof}
By Proposition \ref{doubling} it suffices to prove 
the theorem for $t=1$. We use the simplicial object $\Phi (R\simp^\K)$
for our computations. 
We already know that (\ref{sqodd2}) and (\ref{sqodd3}) hold.
Formula (\ref{sqodd1}) holds for $q=0$ and also for $q=1$ since
$$d_0( \ddr z^{(j)}) = \Sq ^{2j}\ddr y_1 + \binom {n+1} j 
x^j\ddr y_1, \quad 1\leq j \leq n+1.$$
We now prove that (\ref{sqodd1}) also holds for $q=2$. Let 
$v_j$ and $w_j$ denote the following elements for $0\leq j\leq n+1$:
$v_j=\Sq ^{2j} \ddr y_1$ and $w_j=\binom {n+1} j x^j\ddr (y_1)$.
Recall the constructions $m$ and $q$ from section 
\ref{section:normalized}.

{\bf Claim 1}: $m(v_j,v_k)$ is a boundary for $j\neq k$.

We have $d_1v_j=d_1w_j=0$ since $d_1 y_1=0$.
Furthermore, $d_0v_j=0$ and $d_0w_j=0$ since $n$ is odd. 
Thus $v_j$ and $w_j$ are cycles and $v_j+w_j$ is a boundary (actually 
$0$ for $j=0$). Note that
$$m(v_j,v_k)=m(v_j+w_j,v_k)+m(w_j,v_k+w_k)+m(w_j,w_k).$$
The first two terms of the right hand side 
are boundaries by Lemma \ref{lemma:products}. 
We now show that the last term $m(w_j,w_k)$ vanish and Claim 1 follows.

\begin{align*}
m(w_j,w_k) &= s_0(w_j)s_1(w_k)+s_0(w_k)s_1(w_j) \\
&= \binom {n+1} j \binom {n+1} k 
\big( \ddr (x^j y_2) \ddr (x^k y_1) + 
\ddr (x^j y_1) \ddr (x^k y_2) \big)
\end{align*}
so $m(w_j,w_k)=0$ unless $j$ and $k$ are both even. If they are we get
$$m(w_j,w_k) = \binom {n+1} i \binom {n+1} j x^{j+k}
(\ddr y_2 \ddr y_1 + \ddr y_1 \ddr y_2)=0.$$

{\bf Claim 2}: $q(v_j)$ is a cycle for $j>0$ which represents
the same homology class as $\binom {n+1} j x^{2j} \ddr y_1 \ddr y_2$.

By Lemma \ref{lemma:squares} we see that $q(v_j)$ is a cycle.
Note that $q(v_j+w_j)=q(v_j)+q(w_j)+m(v_j,w_j)$. 
By Lemma \ref{lemma:squares} we have that $q(v_j+w_j)$ is a boundary. 
Furthermore, $m(v_j,w_j)$ is a boundary since
$m(v_j,w_j)=m(v_j+w_j,w_j)+m(w_j,w_j)$ and
$m(v_j+w_j,w_j)$ is a boundary by Lemma \ref{lemma:products}
and $m(w_j,w_j)=0$. Thus $q(v_j)+q(w_j)$ is a boundary.
But
\begin{align*}
q(w_j) &= s_0(w_j)s_1(w_j) \\
&= \binom {n+1} j \binom {n+1} j \ddr (x^jy_1) \ddr (x^jy_2) = 
\binom {n+1} j x^{2j} \ddr y_1 \ddr y_2
\end{align*}
so Claim 2 follows. 

We can now prove formula (\ref{sqodd1}) for $q=2$. By the 
Cartan formula we have
\begin{align*}
\Sq^{2i}(\omega_2 ) =& \Sq^{2i}(\ddr y_1 \ddr y_2) \\ 
=& \sum_{j+k=i} \Sq^{2j}(\ddr y_1) \Sq^{2k}(\ddr y_2) \\
=& \sum_{j+k=i, j<k} \big( 
\Sq^{2j}(\ddr y_1) \Sq^{2k}(\ddr y_2)+
\Sq^{2j}(\ddr y_2) \Sq^{2k}(\ddr y_1) \big) \\
&+ \Sq^{i/2} (\ddr y_1) \Sq^{i/2} (\ddr y_2)
\end{align*}
where as usual we put $\Sq^{i/2}=0$ when $i$ is odd. We can 
rewrite this as
$$\Sq ^{2i} (\ddr y_1 \ddr y_2) = 
\sum_{j+k=i,j<k} m(v_j,v_k) + q(v_{i/2}).$$
Now by Claim 1 and Claim 2 see that $\Sq ^{2i}(\ddr y_1 \ddr y_2)$
represents the same homology class as 
$\binom {n+1} {i/2} x^i \ddr y_1 \ddr y_2= 
\binom {2(n+1)} i x^i \omega_2$ for $i$ even and zero
for $i$ odd. Since $\binom {2(n+1)} i =0$ for $i$ odd, the result 
follows.

By the formulas for the product structure and the homotopy operations
we can now prove (\ref{sqodd1}) in general. We start by proving
(\ref{sqodd1}) for $q=2^r$ with $r\geq 1$ by induction on $r$. 
The induction start is the formula for $q=2$, which we just proved. 
Assume that the formula holds for some $q=2^r$.
Put $\gamma_k = \gamma_k (\omega )$ for all $k$.
By Theorem \ref{main1} we have that 
$\gamma_{2q}= \delta_q (\gamma_q )$. 
So we can use the general relations among homotopy operations 
and Steenrod operations 
\cite{G}:
\begin{align*}
\Sq^{2i}\gamma_{2q} &= 
\Sq^{2i} (\delta_q \gamma_q ) =\\
&= \delta_q (\Sq^i \gamma_q )+
\sum_{k<i} \Sq^k (\gamma_q ) 
\Sq^{i-k} (\gamma_q ).
\end{align*}
By induction and by the formula for the product, the last term vanish.
The first term vanish if $i$ is odd. If $i$ is even we get
$$\Sq ^{2i}\gamma_{2q} = 
\binom {q(n+1)} {i/2} \delta_q(x^{i/2} \gamma_q ) =
\binom {q(n+1)} {i/2} x^i \gamma_{2q} .$$
This proves (\ref{sqodd1}) when $q$ is a power of $2$.

Finally, we use the product structure to prove (\ref{sqodd1}) for
all remaining $q$. Assume that it holds for all $q<2^r$ and consider
a $p$ with $2^r<p<2^{r+1}$. We have $p=2^r+s$ where $0<s<2^r$.
Thus $\gamma_p = \gamma_{2^r} \gamma_s$ and
\begin{align*}
\Sq^{2i} (\gamma_p) &= 
\sum_{j=0}^i \Sq ^{2i}(\gamma_{2^r} )\Sq ^{2(i-j)}(\gamma_s ) \\
&=\sum_{j=0}^i \binom {2^r(n+1)} j \binom {s(n+1)} {i-j} x^i\gamma_p
= \binom {p(n+1)} i x^i \gamma_p
\end{align*}   
where we used the combinatorial interpretations of 
binomial coefficients to see the last equality. By induction on $r$
we are done.  
\end{proof}

\begin{theorem}
\label{sqeven}
Let $n\geq 1$, $t\geq 1$ and assume that $n$ is even. 
Then the $\sa$-action on
$H_*(\grtrunc n t x ;\drf)$ is given by the following formulas
where $i,q \geq 0$:
\begin{align}
\label{sqeven1} 
& \Sq^{2^ti} (a_q)= \binom {q(n+1)} i x^i a_q \\
\label{sqeven2}
& \Sq^{2^ti} (b_q)= \binom {q(n+1)+1} i x^i b_q \\
\label{sqeven3}
& \Sq^{2^t i} (x^j) = \binom j i x^{i+j} 
\end{align}
If $k\neq 0$ mod $2^t$ then $\Sq ^k(z)=0$ for all 
$z\in H_*(\grtrunc n t x ;\drf )$.
\end{theorem}

\begin{proof}
It is enough to prove the theorem for $t=1$. 
The proof for the case $n=2$ is given in Lemma 
\ref{sq2viashuffle}. Assume that $n\geq 4$.
We know that (\ref{sqeven3}) holds. Equation (\ref{sqeven1}) and 
(\ref{sqeven2}) are OK for $q=0$ by definition. 

For $q=1$ we use the elements $z^{(i)}$, $1\leq i \leq n$. Since
$$d_0 (\ddr x \ddr z^{(i)}) = \ddr x \Sq ^{2i} \ddr y_1
+ \binom {n+1} i x^i \ddr x \ddr y_1$$
we see that (\ref{sqeven1}) holds for $q=1$.
Equation (\ref{sqeven2}) holds by the following
(where we take $z^{(0)}=0$):

\begin{align*}
& d_0 (x^2\ddr z^{(i-1)} + x \ddr z^{(i)} +z^{(i)} \ddr x) = \\
& x^2 \Sq^{2i-2}(\ddr y_1) + 
\binom {n+1} {i-1} ((i-1) x^i \ddr x y_1+x^{i+1}\ddr y_1) \\
+ & x \Sq^{2i} (\ddr y_1) + 
\binom {n+1} i (ix^i \ddr x y_1 + x^{i+1}\ddr y_1) \\
+ & \Sq^{2i} (y_1) \ddr x +\binom {n+1} i x^i y_1 \ddr x = \\
& \Sq^{2i} (x\ddr y_1 + y_1 \ddr x)+ 
\binom {n+2} i x^i (x\ddr y_1 + y_1\ddr x).
\end{align*}
Note that the assumption, $n$ even, is used for the last
equality.

For $q=2$ we have $\alpha_2=\ddr x \ddr y_1$. By the Cartan formula
we find that
$$\Sq^{2i} (\alpha_2) = \ddr x \sum_{j+k=i,j<k} 
m(\Sq^{2j} \ddr y_1, \Sq^{2k} \ddr y_1)+\ddr xq(\Sq^i \ddr y_1).$$
Define $v_j=\Sq ^{2j} \ddr y_1$ and $w_j= \binom {n+1} j x^j \ddr y_1$.
Formula (\ref{sqeven1}) for $q=2$ now follows by Claim 1 and Claim 2
below.

{\bf Claim 1}: $\ddr x m(v_j,v_k)$ is a boundary when $j\neq k$.

We have that $v_j,w_j\in N_1 \Omega (R\simp )$. Furthermore,
$\ddr xv_j$ and $\ddr xw_j$ are cycles and $v_j+w_j$ is a
boundary. Note that
$$\ddr x m(v_j,v_k)=\ddr x m(v_j+w_j,v_k)
+\ddr m(w_j, v_k+w_k) + \ddr x m(w_j,w_k).$$
The first two terms are boundaries by Lemma \ref{lemma:products} 
and the last term equals zero by direct computation. 
Thus Claim 1 holds.

{\bf Claim 2}: For $i>0$ the element $\ddr x q(\Sq ^{2i}y_1)$ is
a cycle representing the same homology class as 
$\binom {n+1} i x^i \alpha_2$.

We have the following equation:
$$\ddr x q(v_i+w_i)=\ddr x q(v_i)+\ddr x q(w_i) +\ddr x m(v_i,w_i).$$
Since $(v_i+w_i)^2=v_i^2+w_i^2=0$ the left hand side is a boundary
by Lemma \ref{lemma:squares}. The last term on the right hand side
is a boundary since $m(v_i,v_i)=0$ and
$$\ddr x m(v_i,w_i)=\ddr x m(v_i,v_i)+ \ddr x m(v_i, v_i+w_i).$$
Thus $\ddr x q(v_i)$ and $\ddr x q(w_i)$ represent the same
homology class. But $\ddr x q(w_i)=\binom {n+1} i x^i \alpha_2$
so Claim 2 follows. 

We now use the homotopy operations and the product structure
to prove the remaining part of the theorem.
Note that $\Sq^{2i}b_q=C_{q,i}x^i b_q$ and 
$\Sq^{2i}a_q=D_{q,i}x^ia_q$ for some constants 
$C_{q,i}$ and $D_{q,i}$ in $\FF_2$. By the following argument it
suffices to show that the constants are as stated for $1\leq i \leq n-2$.
For $i=0$ the constants are $1$ and for $i=n-1$ they can
can be found from the $i=n-2$ constants 
via the relation $\Sq ^2\Sq ^{2(n-2)}=\Sq ^{2(n-1)}$ which is the
``double'' of $\Sq^1 \Sq^{n-2} = \Sq^{n-1}$. 
For $i\geq n$ there is nothing to show since $x^n=0$.

By the product relation $xa_2=\ddr xb_2$ we see that
\begin{align*}
& C_{2,i} x^{i+1} a_2=C_{2,i} x^i \ddr x b_2=\Sq^{2i} (\ddr x b_2)=
\Sq^{2i} (xa_2)=x^2\Sq^{2i-2}a_2 + x\Sq ^{2i}a_2 \\
&= \binom {2(n+1)} {i-1} x^{i+1}a_2 + \binom {2(n+1)} i x^{i+1}a_2=
\binom {2(n+1)+1} i x^{i+1}a_2.
\end{align*}
Thus formula (\ref{sqeven2}) holds for $q=2$. So we know that 
(\ref{sqeven1}) and (\ref{sqeven2}) holds for $q=0,1,2$.

We now prove that (\ref{sqeven2}) holds when $q=2^r$, $r\geq 1$
by induction on $r$. The initial case $r=2$ holds by the above.
Assume that we have the formula for some $r\geq 2$.
By the general relations \cite{G} we find that
$$\Sq^{2i}(\delta_q(b_q))=\delta_q(\Sq ^ib_q)+
\sum_{t<i} \Sq^t (b_q)\Sq^{i-t}(b_q).$$
The last sum vanish since $b_q^2=0$. The first term vanish when
$i$ is odd. When $i$ is even we use the induction hypothesis:
$$\delta_q(\Sq ^ib_q)=
\binom {q(n+1)+1} {i/2} \delta_q (x^{i/2} b_q) =
\binom {q(n+1)+1} {i/2} x^{i+1}b_{2q}.$$
Thus the following formula holds for all $i$:
$$\Sq^{2i} (\delta_q(b_q))=\binom {2q(n+1)+2} i x^{i+1} b_{2q}.$$
But $\delta_q(b_q)=xb_{2q}$ so this equals
$x\Sq^{2i}(b_{2q})+x^2\Sq^{2i-2}(b_{2q})$ and we have
$$C_{2q,i-1}+C_{2q,i}= \binom {2q(n+1)+2} i.$$ 
Since $C_0=1$ we see that (\ref{sqeven2}) holds for $q=2^{r+1}$.

We now know that (\ref{sqeven2}) holds when $q$ is a power of $2$.
We use the product structure to prove (\ref{sqeven2}) for 
the remaining values of $q$. Assume that the formula holds
for all $q<2^r$ and consider a $p$ with $2^r<p<2^{r+1}$.
Then $p=2^r+s$ where $0<s<2^r$. We have
\begin{align*}
& (C_{p,i}+C_{p,i-1})x^{i+1}b_p =
x\Sq^{2i} b_p+x^2\Sq^{2i-2} b_p = \Sq^{2i} (xb_p)=
\Sq^{2i} (b_{2^r}b_s) \\
&= \sum_{j=0}^i \Sq^{2j} (b_{2^r})\Sq^{2(i-j)}(b_s)
= \sum_{j=0}^i \binom {2^r(n+1)+1} j \binom {s(n+1)+1} {i-j} 
x^ib_{2^r}b_s. 
\end{align*}
By computing the sum of the binomial coefficients we find that
$$D_{p,i}+D_{p,i-1}=\binom {p(n+1)+2} i,$$ 
and since $C_{p,0}=1$ the result follows.

Finally, we finish the proof of (\ref{sqeven1}):
\begin{align*}
& (D_{q,i}+D_{q,i-1})x^{i+1}a_q = x\Sq^{2i} a_q + x^2\Sq^{2i-2}a_q =
\Sq^{2i} (xa_q) = \Sq^{2i} (\ddr x b_q) \\
&= \ddr x \Sq^{2i} b_q = \binom {q(n+1)+1} i \ddr x x^i b_q = 
\binom {q(n+1)+1} i x^{i+1}a_q.
\end{align*}
The result follows since $D_{q,0}=1$.
\end{proof}

The homotopy operations and the product structure 
does not determine the Steenrod algebra action on 
$H_*(\grtrunc n t x;\drf )$ when $n=2$. We prove the
formulas for the action on the chain level instead. 

Let $R\simp^\K \in s\K$ be a simplicial resolution
of $\grtrunc 2 0 x$. As above we have elements
$x, y_1 ,\dots ,y_q\in R_q^\K$ for $q\geq 0$
and $z^{(1)}\in R_2^\K$.

\begin{lemma}
\label{sq2viashuffle}
Put $S\simp = \Omega \Phi R\simp^\K$ 
and let $D$ denote the shuffle map and $\mu$ the product map
as follows:
$$\begin{CD}
S_2\otimes S_{q-2} @>D>> S_q\otimes S_q @>\mu>> S_q.
\end{CD}$$
Define elements $A_q, B_q \in S_2 \otimes S_{q-2}$ for $q\geq 2$ by
$$
A_q= \ddr x\ddr z^{(1)} \otimes \omega_{q-2} \quad , \quad
B_q= \ddr (xz^{(1)}) \otimes \omega_{q-2} +
\sum_{j=1}^{q-2} \ddr x \ddr z^{(1)} \otimes y_j(\omega_{q-2})_j.
$$
Then one has that $\mu D(A_q), \mu D(B_q) \in N_q(S)$ and that
\begin{align*}
& d_0 \mu D(A_q)=
\Sq^2 \alpha_{q-1}+(q-1) x\alpha_{q-1}, \\
& d_0 \mu D(B_q)=
\Sq^2 \beta_{q-1}+q x\beta_{q-1}.
\end{align*}
\end{lemma}

\begin{proof}
By Proposition \ref{EZ} it follows that 
$d_i\mu D(A_q)=0$ for $i>0$ and that
\begin{align*}
d_0\mu D(A_q)&=\mu D(d_0\otimes id+id\otimes d_0)(A_q) \\
&= \mu D(\ddr x \ddr (\Sq^2 (y_1)+xy_1) \otimes \omega_{q-2})+
\mu D(\ddr x\ddr z^{(1)} \otimes x^2\ddr x \omega_{q-3}) 
\end{align*}
Here the last term vanish since $(\ddr x)^2=0$. So by the formula
for the shuffle map we find
\begin{align*}
d_0\mu D(A_q)&= \sum_{i=0}^{q-2} 
s_{q-2} \dots \widehat{s_i} \dots s_0
(\ddr x\ddr (\Sq^2 (y_1)+xy_1))s_i(\omega_{q-2}) \\
&= \sum_{j=1}^{q-1}\ddr x(Sq^2(\ddr y_j)+x\ddr y_j)(\omega_{q-1})_j \\
&= \sum_{j=1}^{q-1} \ddr x\Sq^2 (\ddr y_j)(\omega_{q-1})_j + 
\sum_{j=1}^{q-1}x\ddr x \omega_{q-1} \\
&= \Sq^2 \alpha_{q-1} + (q-1) x\alpha_{q-1}.
\end{align*}

We now check that $d_i\mu D(B_q)=0$ for $i>0$. Since 
$d_t(z^{(1)})=d_t(\omega_{q-2})=0$ for $t>0$ we
get the following by Proposition \ref{EZ}:
$$d_i\mu D(B_q) = \sum_{j=1}^{q-2} \sum_{t=1}^i
\mu G_t(\ddr x \ddr z^{(1)} \otimes d_t(y_j(\omega_{q-2})_j)).$$
When $1\leq t \leq q-3$ one has that
$d_t(y_j(\omega_{q-2})_j)= y_j \omega_{q-3}$ for $j=t,t+1$ 
and zero otherwise. When $t=q-2$ we always get zero. 
So $d_i\mu D(B_q)=0$ for $i>0$ as stated.

Finally, we show that 
$d_0\mu D(B_q)=\Sq^2 \beta_{q-1}+qx\beta_{q-1}$. The right
hand side gives the following:
\begin{align}
\label{resultofsum}
& \Sq^2 \beta_{q-1}+qx\beta_{q-1} =  
\nonumber \\
& x\Sq^2(\omega_{q-1} )+x^2\omega_{q-1} +
\ddr x \sum_{j=1}^{q-1} \Sq^2 (y_j)(\omega_{q-1})_j +
\ddr x \sum_{j=1}^{q-1} y_j \Sq^2 ((\omega_{q-1})_j) 
\nonumber \\
&+ qx^2\omega_{q-1} + 
qx\ddr x \sum_{j=1}^{q-1} y_j(\omega_{q-1} )_j = 
\nonumber \\
& x\Sq^2(\omega_{q-1} )+(q+1)x^2\omega_{q-1}
+qx\ddr x \sum_{j=1}^{q-1} y_j(\omega_{q-1} )_j
+\ddr x \sum_{j=1}^{q-1} \Sq^2 (y_j)(\omega_{q-1})_j
\nonumber \\
&+\ddr x \sum_{j=1}^{q-1} y_j
(\sum_{k=1}^{j-1} \Sq^2 (\ddr y_k)(\omega_{q-1})_{k,j} + 
\sum_{k=j+1}^{q-1} \Sq^2 (\ddr y_k)(\omega_{q-1})_{j,k})
\end{align}

For the left hand side, we use that 
$d_0\mu D=\mu D(d_0\otimes id+id\otimes d_0)$ by Proposition
\ref{EZ}. Since $(\ddr x)^2=0$ we find that
$$\mu D(id\otimes d_0)(B_q)=
\mu D(\ddr (xz^{(1)})\otimes x^2\ddr x\omega_{q-3} +
\ddr x \ddr z^{(1)} \otimes x^3 \omega_{q-3})=0$$
and hence
\begin{align}
\label{twoterms}
& d_0\mu D(B_q)=\mu D(d_0\otimes id)(B_q) \nonumber \\
&= \mu D(\ddr (xd_0z^{(1)})\otimes \omega_{q-2})
+\ddr x \sum_{j=1}^{q-2} 
\mu D(\ddr (d_0z^{(1)})\otimes y_j(\omega_{q-2})_j). 
\end{align}
The first term on the right hand side of (\ref{twoterms}) 
gives the following:
\begin{align}
\label{firstterm}
&\mu D(\ddr (xd_0z^{(1)}) \otimes \omega_{q-2} ) = 
\sum_{j=1}^{q-1} \ddr (x\Sq^2 (y_j)+x^2y_j) 
(\omega_{q-1})_j = \nonumber \\
& \ddr x \sum_{j=1}^{q-1} \Sq^2 (y_j)(\omega_{q-1} )_j
+ x\Sq^2 (\omega_{q-1} ) + (q-1)x^2\omega_{q-1}.
\end{align}
By the formula for the shuffle map $D$ the second term on the 
right hand side of (\ref{twoterms}) equals
\begin{align*}
& \ddr x \sum_{j=1}^{q-2} \sum_{k=1}^{q-1}
(\Sq^2 (\ddr y_k)+x\ddr y_k)s_{k-1}( y_j(\omega_{q-2} )_j) =
\\
& \ddr x \sum_{j=1}^{q-2} \sum_{k=1}^{j}
(\Sq^2 (\ddr y_k)+x\ddr y_k)y_{j+1}(\omega_{q-1} )_{k,j+1} 
\\
+& \ddr x \sum_{j=1}^{q-2} \sum_{k=j+1}^{q-1}
(\Sq^2 (\ddr y_k)+x\ddr y_k)y_{j}(\omega_{q-1} )_{j,k}.
\end{align*}
After performing the substitution $j^\prime = j+1$
in the first sum on the right hand side we see that this equals
\begin{align*}
& x\ddr x \sum_{j=2}^{q-1} (j-1)y_j(\omega_{q-1} )_j 
+ x\ddr x \sum_{j=1}^{q-2} (q-1-j)y_j(\omega_{q-1} )_j 
\\
+& \ddr x \sum_{j=2}^{q-1} \sum_{k=1}^{j-1}
y_j \Sq^2 (\ddr y_k)(\omega_{q-1})_{k,j} 
+\ddr x \sum_{j=1}^{q-2} \sum_{k=j+1}^{q-1}
y_j \Sq^2 (\ddr y_k)(\omega_{q-1})_{j,k}.
\end{align*}
By a final reduction we conclude that the second term 
on the right hand side of (\ref{twoterms}) equals
\begin{align}
\label{secondterm}
& (q-2)x\ddr x \sum_{j=1}^{q-1} y_j (\omega_{q-1} )_j
\nonumber \\
+&\ddr x \sum_{j=1}^{q-1} y_j (
\sum_{k=1}^{j-1} \Sq^2 (\ddr y_k)(\omega_{q-1} )_{k,j} 
+ \sum_{k=j+1}^{q-1} \Sq^2 (\ddr y_k)(\omega_{q-1} )_{j,k})
\end{align}
Since the sum of (\ref{firstterm}) and  (\ref{secondterm}) 
equals (\ref{resultofsum}) we are done.
\end{proof}

When $n=1$, $|x|\geq 1$ we have that 
$\trunc n x = \Lambda (x)$ which is an object in 
$\K$ with trivial $\sa$-action. 

\begin{proposition}
\label{spherecase}
The $\sa$-action on $H_*(\Lambda (x);\drf )$ is trivial
when $|x|\geq 2$.
\end{proposition}

\begin{proof}
Put $m=|x|$. We have computed 
$H_*(\Lambda (x) ;\drf ) \cong 
\Lambda (x) \otimes \Lambda (\ddr x) \otimes \Gamma [\omega ]$
where $||x||=(m,0)$, $||\ddr x|| = (m-1,0)$ and
$||\gamma_q (\omega )|| = (-q,q(2m-1))$.

Because of the bigrading the only possible non zero 
squares on $\gamma_q (\omega )$ are
$\Sq^{m-1} \gamma_q (\omega )$, $\Sq^{m} \gamma_q (\omega )$,
$\Sq^{2m-1} \gamma_q (\omega )$ which might equal
$\ddr x \gamma_q (\omega )$, $x \gamma_q (\omega )$,
$x\ddr x \gamma_q (\omega )$ respectively.
Recall that $Sq^r$ is decomposable unless $r=2^s$, $s\geq 0$.
So the result follows for all $m$ except for
$m=2^k$ and $m=2^k+1$ with $k\geq 1$. When $m=2^k$ the
result is a special case of Theorem \ref{sqodd}.

Assume that $m=2^k+1$, $k\geq 1$. Here we must  
show that $\Sq^{m-1} \gamma_q (\omega )=0$. By the 
homotopy operations and algebra structure we see that
it suffices to show that 
$\Sq^{m-1} \gamma_1 (\omega )=0$ and
$\Sq^{m-1} \gamma_2 (\omega )=0$.

Let $R\simp^\K$ be a simplicial resolution of $\Lambda (x)$ 
in $s\K$ such that $R_0^\K$ is the free object in $\K$
on the generator $x$ (thus $R_0^\K \cong H^*(K(\FF_2 ,m))$).
We have a weak equivalence of simplicial algebras
$R\simp \to R\simp ^\K$. Let $y_1\in R\simp^\K$ denote the 
image of $y_1 \in R\simp$. We have $d_0y_1=x^2$ and
$d_1y_1=0$. Since $\pi_0(R\simp^\K)=\Lambda (x)$ 
there are elements $w^{(i)} \in R_1^\K$ for $1\leq i \leq m-1$
such that $d_0w^{(i)}=\Sq^i x$ and $d_1w^{(i)}=0$.

But $d_0(\Sq^{2i} y_1+(w^{(i)})^2)=\Sq^{2i}(x^2)+(\Sq^i x)^2=0$
and $d_0(\Sq^{2i+1} y_1)=\Sq^{2i+1} (x^2)=0$ so there 
exist elements $z^{(t)} \in \R^\K _2$ for $1\leq t \leq 2m-1$ 
with $d_jz^{(t)}=0$ for $j=1,2$ and 
$d_0z^{(2i)}=\Sq^{2i} y_1 +(w^{(i)})^2$,
$d_0z^{(2i+1)}=\Sq^{2i+1} y_1$.

In $R\simp^\K$ we have $\ddr z^{(t)}$ with 
$d_j \ddr z^{(t)}=0$ and $d_0 \ddr z^{(t)} = \Sq^t \ddr y_1$.
Thus $\Sq^{m-1} \gamma_1 (\omega )=0$. For $\gamma_2 (\omega)$
note that
$$\Sq^{2^k}(\ddr y_1 \ddr y_2) = \sum_{t=0}^{2^{k-1}-1} 
m(d_0\ddr z^{(t)},d_0\ddr z^{(m-t)}) 
+ q(d_0 \ddr z^{(2^{k-1})})$$
so the result follows by Lemma \ref{lemma:products} and 
Lemma \ref{lemma:squares}.
\end{proof}

\section{The results}

Recall that if $X$ is a space with 
$H^*X=\FF_2 [x]/(x^{n+1})$ where $n\geq 2$ then
$|x|=2^k$ for $k=0,1,2$ or $|x|=8$ and $n=2$ \cite{A}.
Examples are the projective spaces $\rp^n$, $\cp^n$, $\hp^n$ 
and the Cayley projective plane $\cayley$, but there are other 
homotopy types of spaces with these cohomology rings.
For $n=1$ there is a space $X$ with $H^*X=\Lambda (x)$ 
for any positive value of $|x|$. The spheres $S^m$ for $m\geq 1$
are examples of such spaces.

When $n$ is odd one cannot determine the operation $\Sq ^1$ by 
our spectral sequence. For the projective spaces and spheres this
information follows from the results in \cite{Z}.

\begin{proposition}
\label{Ziller}
Consider the operation 
$\Sq^1 : H^t(\Lambda X) \to H^{t+1}(\Lambda X)$. 

1) $X=\cp^n$, $n\geq 1$: If $n\neq 1$ mod $4$ then all $\Sq^1=0$.
If $n=1$ mod $4$ then the only non vanishing $\Sq^1$'s appear when 
$t\in \{ 2nk|k\geq 1 \}$ and here one has
$\Sq ^1=(0, id): \FF_2 \oplus \FF_2 \to \FF_2 \oplus \FF_2$.

2) $X=\hp^n$, $n\geq 1$: If $n\neq 1$ mod $4$ then all $\Sq^1=0$.
If $n=1$ mod $4$ then the only non vanishing
$\Sq^1$'s appear when 
$t\in \{ (4n+2)k|k\geq 1 \}$ 
and here $\Sq ^1=id: \FF_2 \to \FF_2$.

3) $X=S^m$, $m\geq 3$: If $m$ is odd then all $\Sq^1=0$. 
If $m$ is even then the only non vanishing $\Sq^1$'s appear when 
$t\in \{2k(m-1) |k\geq 1 \}$ and here
$\Sq ^1=id: \FF_2 \to \FF_2$. 
\end{proposition}

\begin{remark} 
In 3) we do not mention the case $m=2$. This case is however 
included in 1) since $S^2=\cp^1$.
\end{remark}

\begin{proof}
The homology groups of $\Lambda S^m$, $\Lambda \cp^n$ and
$\Lambda \hp^n$ with $\ZZ$-coefficients are listed in 
\cite{Z} on page 21.
By universal coefficients one finds the cohomology
groups with coefficients in $\ZZ /2 $ and $\ZZ /4$.
The result follows by the long exact coefficient sequence
associated with $0 \to \ZZ /2 \to \ZZ /4 \to \ZZ/2 \to 0$.
\end{proof}

\begin{theorem}
\label{main2}
Assume that $X$ is a 1-connected space with cohomology algebra
$H^*X=\FF_2 [x]/(x^{n+1})$ where $n$ is a positive integer. 
Put $m=|x|$ and recall that when $n=1$, $n=2$, $n\geq 3$ then 
$m\in \ZZ_+$, $m\in \{ 2,4,8 \}$, $m\in \{ 2,4 \}$ respectively. 

1) Assume that $n$ is even. Then there is an isomorphism
$$H^*(\Lambda X) \cong \FF_2 \oplus 
\FF_2 [x]/(x^n) <a_q,b_q|q\geq 0>$$
where $|x|=m$, $|a_q|=qm(n+1)-2q+m-1$ and $|b_q|=qm(n+1)-2q+m$.
The algebra structure is given by the relations
$$a_pa_q=0,\quad b_pb_q = \binom {p+q} p xb_{p+q}, \quad
a_pb_q = \binom {p+q} p xa_{p+q}.$$ 
The Steenrod algebra structure is given as follows where $i,j,q\geq 0$:
\begin{align*}
& \Sq^{mi} (x^ja_q) = \binom {q(n+1)+j} i x^{i+j} a_q, \\
& \Sq^{mi} (x^jb_q) = \binom {q(n+1)+1+j} i x^{i+j} b_q
\end{align*}
and $\Sq^t z = 0$ for all $z\in H^*(\Lambda X)$ when 
$t\neq 0$ mod $m$.

2) Assume instead that $n$ is odd.
Then there is an isomorphism of algebras
$$H^*(\Lambda X) \cong 
\FF_2 [x]/(x^{n+1}) \otimes \Lambda (\ddr x) \otimes \Gamma [\omega ]$$
where $|x|=m$, $|\ddr x|=m-1$, $|\gamma_q (\omega )|=qm(n+1)-2q$.
The Steenrod algebra action satisfies the following for $i,j,q\geq 0$
and $\epsilon = 0,1$:
$$
\Sq^{mi} (x^j(\ddr x)^\epsilon \gamma_q (\omega )) = 
\binom {q(n+1)+j} i x^{i+j}(\ddr x)^\epsilon \gamma_q (\omega ) 
$$
If $t\neq 0$ mod $m$ then
$\Sq^t z =0$ for all $z\in H^*(\Lambda X)$ unless   
$z= \gamma_q (\omega )$ for a $q\geq 1$ and $t=1$. Here 
$\Sq^1 \gamma_q (\omega )=c_qx^n\ddr x \gamma_{q-1} (\omega )$
for a constant $c_q \in \FF_2$. 
If $X=\cp^n$ or $X=\hp^n$ then $c_q=1$ for $n=1$ mod $4$ and 
$c_q=0$ for $n=3$ mod $4$, for all $q\geq 1$.
If $X=S^m$, $m\geq 2$ then $c_q=0$ for $m$ odd and 
$c_q=1$ for $m$ is even, for all $q\geq 1$.
\end{theorem}

\begin{proof}
The distribution of zero's in the $E_2$ therm shows that 
all differentials are trivial. So the spectral sequence
collapse and $E_2=E_\infty$. 
When $n$ is even, there is
a unique representative in $H^*(\Lambda X)$ for each
element in $E_\infty$ and the result follows easily.

Assume that $n$ is odd. We have a filtration of 
$H^*=H^*(\Lambda X)$ as follows:
$$H^*\supseteq \dots \supseteq F^{-2}H^* \supseteq F^{-1}H^*
\supseteq F^{0}H^* \supseteq F^{1}H^*=0$$
with $F^{-p}H^{q-p}/F^{-p+1}H^{q-p} \cong E_\infty^{-p,q}$.
Define an algebra map 
$$\phi : 
\trunc n x \otimes \Lambda (\ddr x) \otimes \Gamma [\omega ]
\to H^*$$
by $x\mapsto x$, $\ddr x \mapsto \ddr x$ and
$\gamma_i (\omega ) \mapsto \gamma_i^\prime $
where $x$ and $\ddr x$ are the unique representatives
for $x$ and $\ddr x$ in $E_\infty$ and 
$\gamma_i^\prime$ is one of the two representatives
for $\gamma_i (\omega )$ in $E_\infty$. We must check 
that $\phi$ is well defined.

Since $\Gamma [\omega ]\cong 
\Lambda (\gamma_q (\omega )| q=2^j, j\geq 0)$ it suffices
to see that $(\gamma_q^\prime )^2=0$ when $q$ is a power of
two. From the $E_\infty$-term we find that
\begin{align*}
(\gamma_q^\prime )^2 =\Sq^{q(m(n+1)-2)} (\gamma_q^\prime ) &=
\binom {mq(n+1)} {mq(n+1)-2q} 
x^{q(n+1)-2q/m} \gamma_q^\prime \\
&= \binom {m(n+1)q} {2q} 
x^{q(n+1)-2q/m} \gamma_q^\prime =0
\end{align*} 
So $\phi$ is well defined. We can now define a filtration 
on its domain space such that $\phi$ becomes a map of 
filtered rings. Since the associated map of graded objects 
is an isomorphism, $\phi$ is an isomorphism.

It follows directly from the $E_\infty$ term and Proposition
\ref{Ziller} that the Steenrod algebra action is as stated.
\end{proof}

\section{A Thom space interpretation}

Let $M^d$ be a closed, connected and oriented smooth manifold
without boundary. We denote the dual cohomology class of the 
fundamental class $[M] \in H_d(M;\ZZ )$ by $u_0$. 
Let $\tau \to M$ be the tangent bundle with Thom space 
$\thom \tau$ and Thom class $u_\tau \in H^d(M;\ZZ )$. 
The zero section defines a map $s_0:M\to \thom \tau$.
Recall that $s_0^*: H^d(\thom \tau ;\ZZ ) \to H^d(M;\ZZ )$
is given by multiplication with the Euler characteristics
$$s_0^*(u_\tau ) = \chi (M) u_0.$$ 

Consider the $q$ fold Whitney sum 
$q\tau= \tau \oplus \dots \oplus \tau$. 
We view $q\tau$ as a sub bundle of $(q+1)\tau$ by 
the identity on the first $q$ summands and zero on the 
last summand. Thus we have an
inclusion $s_q: \thom {q\tau } \to \thom {(q+1)\tau }$.

\begin{lemma}
\label{Euler}
The map $s_q^*:H^{(q+1)d}(\thom {(q+1)\tau } ;\ZZ ) \to
H^{(q+1)d}(\thom {q\tau } ;\ZZ )$
is given by multiplication with the Euler characteristics
$$s_q^*(u_{(q+1)\tau }) = \chi (M) u_{q\tau }u_0.$$ 
\end{lemma}

\begin{proof}
The Whitney sum $(q+1)\tau$ is the pullback of the product
bundle $\tau^{\times (q+1)} \to M^{\times (q+1)}$
along the diagonal $\Delta_{q+1} :M \to M^{\times (q+1)}$.

We write $\epsilon^0 \to M$ for the trivial 0-dimensional
vector bundle over $M$. This bundle simply means the 
identity map $M\to M$ and its Thom space is $M$ with a
disjoint base point $\thom {\epsilon^0} = M_+$.

Form the pullback of $\tau^{\times q} \to M^{\times q}$
along a composite map as follows:
$$
\begin{CD}
q\tau @>>> \tau^{\times q} 
@>>> \tau^{\times q} \times \epsilon^0
@>>> \tau^{\times q} \\
@VVV @VVV @VVV @VVV \\
M @>\Delta_q>> M^{\times q} 
@>id \times \Delta_1>> M^{\times (q+1)} 
@>pr>> M^{\times q}
\end{CD}
$$
where the projection $pr$ is on the first $q$ factors.
The pulback of $\tau^{\times q}$ along $pr$ equals
$\tau^{\times q} \times \epsilon^0$. The further pulback
of this vector bundle is as stated in the diagram since
$pr \circ (id\times \Delta_1)=id$ 
and $(id\times \Delta_1)\circ \Delta_q = \Delta_{q+1}$.

The vector bundle $\tau^{\times q} \times \epsilon^0$ is
a sub bundle of $\tau^{\times (q+1)}$ by the map
$id \times s_0$. Furthermore 
$\tau^{\times q} \times \epsilon^0 \subset
\tau^{\times (q+1)}$
pulls back to $q\tau \subset (q+1)\tau$ along
$\Delta^{q+1}$. So we get a commutative diagram
as follows (with $\ZZ$-coefficients):
$$
\begin{CD}
H^*({\thom \tau }^{\wedge (q+1)}) @>\cong>> 
H^*(\thom {\tau^{\times (q+1)} }) @>>>   
H^*(\thom {(q+1)\tau } ) \\
@VVV @V(id\times s_0)^*VV @VVV \\
H^*({\thom \tau }^{\wedge q} \wedge M_+ ) @>\cong>>
H^*(\thom {\tau^{\times q} \times \epsilon^0 };\ZZ ) @>>>
H^*(\thom {q\tau } )
\end{CD}
$$
The Thom class for $\tau^{\times (q+1)}$ 
is $u_\tau ^{\wedge (q+1)}$ and it maps to
the Thom class $u_{(q+1)\tau}$ by the
upper horizontal map. But 
$(id\times s_0)^* (u_\tau ^{\wedge (q+1)})=
u_\tau^{\wedge q}\wedge \chi (M) u_0$ so
the result follows by the Thom isomorphism
for the lower horizontal map.
\end{proof}

\begin{definition}
\label{cofdef}
Let $C_q(M)$ denote the cofiber of $s_q$:
$$\thom {q\tau } \to \thom {(q+1)\tau } \to C_q(M)$$
\end{definition}

\begin{proposition}
\label{cofcom}
Let $M^d$ be a compact, connected smooth 
manifold without boundary and of dimension $d\geq 2$. 
Assume that $H_1(M;\ZZ)=0$ and that $H_{d-1}(M;\ZZ )$ is 
torsion free. Assume also that the Euler characteristics
$\chi (M)\neq 0$. Put $T_q = \thom {q\tau }$, $C_q= C_q(M)$ and
write $\phi$ for the Thom isomorphism and $\partial$ for
the connecting homomorphism. \\
1) $\ZZ$-coefficients: There are isomorphisms
\begin{align*}
& \rhz k {C_q} \xrightarrow[\cong]{\partial_*}
\rhz {k-1} {T_q} \xrightarrow[\cong]{\phi_*}
\hz {k-1-qd} M  
&,& \quad k<(q+1)d \\
& \rhz k {C_q} \cong \ZZ / \chi (M)
&,& \quad  k=(q+1)d \\
& \rhz k {C_q} \xleftarrow[\cong]{} 
\rhz k {T_{q+1}} \xrightarrow[\cong]{\phi_*}
\hz {k-(q+1)d} M 
&,& \quad k>(q+1)d.
\end{align*}
2) $\FF_2$-coefficients: There are isomorphisms
\begin{align*}
& \rcoh k {C_q} \xleftarrow[\cong]{\partial^*}
\rcoh {k-1} {T_q} \xleftarrow[\cong]{\phi}
\coh {k-1-qd} M 
&,& \quad k < (q+1)d\\
& \rcoh k {C_q} \xrightarrow[\cong]{} 
\rcoh k {T_{q+1}} \xleftarrow[\cong]{\phi}
\coh {k-(q+1)d} M 
 &,& \quad k> (q+1)d+1.
\end{align*}
If $\chi (M)$ is odd then
$\rcoh {(q+1)d} {C_q} = \rcoh {(q+1)d+1} {C_q} =0$.
If $\chi (M)$ is even then
$$\rcoh {(q+1)d} {C_q} \xrightarrow[\cong ]{} 
\rcoh {(q+1)d} {T_{q+1}} \quad , \quad 
\rcoh {(q+1)d+1} {C_q} \xleftarrow[\cong ]{\partial^*}
\rcoh {(q+1)d} {T_q}.$$
The operation 
$\Sq^1 :\rcoh {(q+1)d} {C_q} \to 
\rcoh {(q+1)d+1} {C_q}$
is the identity map $id:\FF_2 \to \FF_2$ when $\chi (M)=2$ mod $4$ and
the trivial map $0:\FF_2 \to \FF_2$ when $\chi (M)=0$ mod $4$. 
\end{proposition} 

\begin{proof}
By universal coefficients $H^1(M;\FF_2)=0$ so $M$ is oriented.
Poincar\'{e} duality gives that $H^{d-1}(M;\ZZ )=0$. So we have
that $H_{d-2}(M;\ZZ )$ is torsion free and 
$H_{d-1}(M;\ZZ )=0$ by universal coefficients.
It follows that $H^{d-1}(M;\FF_2 )=0$. 
 
The long exact (co)homology sequence and the Thom
isomorphism gives the results except for the statement
regarding $\Sq^1$. By universal coefficients one
can determine the cohomology groups with $\ZZ /4$-coefficients.
The long exact coefficient sequence associated with
$0\to \ZZ /2 \to \ZZ /4 \to \ZZ /2 \to 0$
then determines $\Sq^1$ since it agrees with the 
Bockstein homomorphism.
\end{proof}

\begin{theorem}
\label{mainresult}
Let $M$ be a complex or quaternion projective space
$\cp^n$, $\hp^n$, $n\geq 1$ the Cayley projective plane 
$\cayley$ or a sphere $S^m$, $m\geq 2$. 
Let $r(M)$ denote the degree of the algebra generator for 
$H^*(M;\FF_2 )$. Define the pointed space $\lm (M)$ by
$$\lm (M) = M_+ \vee 
\bigvee_{q\geq 0} \Sigma^{(r(M)-2)(q+1)} C_q(M).$$
Then the suspension spectra
$\Sigma^\infty (\Lambda M)_+$ and $\Sigma^\infty \lm (M)$ 
have isomorphic homology groups with $\ZZ$-coefficients.
Furthermore the cohomology groups with $\FF_2$-coefficients
are isomorphic modules over the Steenrod algebra.
For the spheres we have a homotopy equivalence
$\Sigma^\infty (\Lambda S^m)_+ \simeq \Sigma^\infty CT(S^m)$.
\end{theorem}

\begin{proof}
We first prove the statement regarding the projective spaces and
$\ZZ$-coefficients. Let $\num K = \num C, \num H, \numcayley$ and put 
$r=r(M)=\dim_{\num R} \num K$.

We use Proposition \ref{cofcom} to compute the stable 
homology of $\lm (M)$. The non vanishing groups 
$\rhz k {C_q(\kp^n )}$ equal $\ZZ$ for
$$k\in \{ qd+1+jr, (q+1)d+(j+1)r | 0 \leq j \leq n-1 \} ,$$ 
and $\ZZ /(n+1)$ for $k=(q+1)d$. Thus 
$\rhz k {\Sigma^{(r-2)(q+1)} C_q(\kp^n )}$ 
equals $\ZZ$ for $k \in A(r,q)\cup B(r,q)$ 
and $\ZZ / (n+1)$ for $k=\big( (n+1)r-2 \big) (q+1)$ 
and $0$ otherwise where
\begin{align*}
& A(r,q) = \{ \big( (n+1)r-2\big) q
+rj+r-1 |0\leq j \leq n-1 \} , \\
& B(r,q) = \{ \big( (n+1)r-2\big) (q+1)
+r(j+1) |0\leq j \leq n-1 \} .
\end{align*} 
Finally, $\hz k {\kp^n }=\ZZ$ for $k\in C(r)$ 
where $C(r)=\{ ri|0\leq i \leq n \}$ and zero otherwise.

For $\num K= \num C$, $r=2$ we have that
$C(2)\cup A(2,0)= \{ 0 , \dots , 2n \}$ and that 
$B(2,q-1)\cup A(2,q)= \{ 2nq+1 , \dots ,2n(q+1) \}$ for
$q\geq 1$. So the group $\rhz k {\lm (\cp^n )}$ equals 
$\ZZ \oplus \ZZ /(n+1)$ for $k=2nm$, $m\geq 1$ and
$\ZZ$ otherwise.
For $\hz k {\Lambda \cp^n}$ we have the same by the
table in \cite{Z}.

For $\num K= \num H$, $r=4$ the group $\rhz k {\lm (\hp^n )}$ 
equals $\ZZ$ when $k$ belongs to one of the sets
\begin{align*}
& C(4)\cup \cup_{q\geq 0} B(4,q) = 
\{ 0\} \cup \{ 2m(2n+1)+4l| m\geq 1 , 1\leq l \leq n \} ,\\
& \cup_{q\geq 0} A(4,q) =
\{ 2m(2n+1)+4l-4n+1|m\geq 1, 0\leq l \leq n-1 \} ,
\end{align*}
and $\ZZ / (n+1)$ when $k\in \{ 2m(2n+1) | m\geq 1 \}$. 
Otherwise the group vanish. For $\hz k {\Lambda \hp^n }$ we have 
the same by \cite{Z}.

For $\num K = \numcayley$, $r=8$, $n=2$ the group
$\rhz k {\lm (\cayley )}$ equals $\ZZ$ for  
$k\in \{ 0,8,16 \} \cup \{ 22m-15, 22m-7, 22m+8, 22m+16|m\geq 1 \}$
and a $\ZZ /3$ for $k\in \{ 22m|m \geq 1\}$.
The list in \cite{Z} differs here by having zero's in 
degrees $0,8,16$ and $22m-15$, $m\geq 1$ instead of $\ZZ$'s.

This is because there is an error in Ziller's final table of
$H_*(\Lambda \cayley ;\ZZ)$. His  
result is the same as ours when one corrects this error as we 
will now see. By \cite{Z} Theorem 8 one has
$$H_*(\Lambda \cayley ;\ZZ ) \cong 
\bigoplus_{B} H_{*-\lambda (B)}(B; \ZZ ),$$
where the sum is over all critical sub manifolds $B\subset \Lambda \cayley$
and $\lambda (M)$ is the index of $M$. All the critical sub manifolds
are diffeomorphic to the sphere bundle $T_1M=S(\tau )$. By the table 
on page 20 of \cite{Z} one sees that
$H_k(T_1\cayley ;\ZZ )= \ZZ$ for $k=0,8,23,31$ and 
$H_k(T_1\cayley ;\ZZ )=\ZZ /3$ for $k=15$ and zero otherwise.
The index can be found in the table on page 11: 
$\lambda (c^m)=22m-15$.

So Ziller shows that
$$H_k(\Lambda \cayley ; \ZZ ) \cong H_k(\cayley ; \ZZ) \oplus 
\bigoplus_{m\geq 1} H_{k-22m+15} (T_1 \cayley ;\ZZ ).$$
Thus one gets a $\ZZ$ when 
$k\in \{ 0,8,16 \}$ or $k=22m-15+t$ for $t\in \{ 0,8,23,31 \}$,
$m\geq 1$ and one gets a $\ZZ /3$ when $k=22m$, $m\geq 1$. 
This is the same as our result.

Now we prove the statements regarding $\FF_2$-coefficients
for the projective spaces. We use Proposition \ref{cofcom} to
find the $\sa$-action on $\rcoh * {\lm (\kp^n )}$.
We have that 
$H^*(\kp^n )= \FF_2 [x]/(x^{n+1})$ where $|x|=r$.
By \cite{MS} page 133, Wu's theorem implies that the total 
Stiefel-Whitney class of $\tau$ is
$w(\tau )=(1+x)^{n+1}$ and hence
$$w(q\tau )=(1+x)^{q(n+1)} = 
\sum_{i=0}^\infty \binom {q(n+1)} i x^i.$$ 
We have a Thom isomorphism 
$\phi : H^*(\kp^n ) \to \rcoh {*+rnq} {\thom {q\tau }}$
and the action on the fundamental class satisfies 
$\Sq^k u_{q\tau} = \phi (w_k (q\tau))$ such that
\begin{equation}
\label{sqfund}
Sq^{ri} u_q = \binom {q(n+1)} i x^i u_q \quad , \quad
\Sq^{k} u_q = 0 \text{ if } k \neq 0 \text{ mod } r
\end{equation}

Assume that $n$ is even. Then $\chi (\kp^n )=n+1$ is odd.
Proposition \ref{cofcom} shows that 
$\rcoh * {\Sigma^{(r-2)(q+1)} C_q(\kp^n )}$
has $\FF_2$-basis $\{ a_q^j , b_{q+1}^j |0\leq j <n \}$ with
\begin{align*}
& a_q^j=\sigma^{(r-2)(q+1)} \partial^* (x^j u_{q\tau })
\quad , \quad  
b_{q+1}^j= \sigma^{(r-2)(q+1)} (Q^*)^{-1} (x^{j+1}u_{(q+1)\tau }), \\
& |a_q^j|=qr(n+1)-2q+r-1+rj, \\
& |b_{q+1}^j|=(q+1)r(n+1)-2(q+1)+r+rj 
\end{align*}
where $Q:\thom {(q+1)\tau} \to C_q (\kp^n )$ denotes the 
quotient map.
By (\ref{sqfund}) and the Cartan formula we have
\begin{align*}
& \Sq^{ri} (a_q^j) = \binom {q(n+1)+j} i a_q^{i+j} , \\
& \Sq^{ri} (b_{q+1}^j) = \binom {(q+1)(n+1)+1+j} i b_{q+1}^{i+j}
\end{align*}
where by definition $a_q^t=b_{q+1}^t=0$ if $t\geq n$. The operations
$\Sq^k$ for $k\neq 0$ mod $r$ sends the basis elements to zero.
The reason why there is no therms of the form $a_q^j$ in 
the formula for $\Sq^t(b_{q+1}^j)$ is that the degrees of
such terms are to small.
The result follows by comparing with Theorem \ref{main2}.

Assume that $n$ is odd. Then $\chi (\kp^n )=n+1$ is even.
Proposition \ref{cofcom} shows that
$\rcoh * {\Sigma^{(r-2)(q+1)} C_q(\kp^n )}$
has $\FF_2$-basis $\{ c_q^j , d_{q+1}^j |0\leq j \leq n \}$ with
\begin{align*}
& c_q^j=\sigma^{(r-2)(q+1)} \partial (x^j u_{q\tau })
\quad , \quad  
d_{q+1}^j= \sigma^{(r-2)(q+1)} Q^{-1} (x^ju_{(q+1)\tau }), \\
& |c_q^j|=qr(n+1)-2q+1+rj, \\
& |d_{q+1}^j|=(q+1)r(n+1)-2(q+1)+rj 
\end{align*}
By (\ref{sqfund}) and the Cartan formula we have
\begin{align*}
& \Sq^{ri} (c_q^j) = \binom {q(n+1)+j} i c_q^{i+j} , \\
& \Sq^{ri} (b_{q+1}^j) = \binom {(q+1)(n+1)+j} i b_{q+1}^{i+j}
\end{align*}
where by definition $a_q^t=b_{q+1}^t=0$ if $t>n$. The operations
$\Sq^k$ for $k\neq 0$ mod $r$ vanish on the basis elements except for
$\Sq^1 d_{q+1}^0$ which equals $c_q^n$ when $n=1$ mod $4$ and
zero when $n=3$ mod 4. Note that we cannot hit any of the classes 
$b_{q+1}^j$ by a Steenrod operation on one of the classes
$a_q^j$ since $a_q^j$ is in the image of $\partial$.
The result follows by comparing with Theorem \ref{main2}. 
The class $c_q^j$ corresponds to $x^j\ddr x \gamma_q (\omega )$
and $d_{q+1}^j$ corresponds to $x^j\gamma_{q+1} (\omega )$.

For the spheres we can prove directly that there is a 
homotopy equivalence. We first identify the stable homotopy type of
$C_q(S^m)$. We have $S^m\subseteq {\num R}^{m+1}$ with
trivial one dimensional normal bundle $\epsilon^1 \to S^m$ and
$\tau \oplus \epsilon^1 \cong \epsilon^{m+1}$. Thus
\begin{align*}
\Sigma^{q+t} \thom {q\tau} \cong &
\thom {(q+t)\epsilon^1 \oplus q\tau } \cong
\thom {\epsilon^{q(m+1)+t}} \cong 
\Sigma^{q(m+1)+t}S^m_+ \cong \\
& S^{q(m+1)+t}\vee S^{q(m+1)+m+t}.
\end{align*}
By taking the $(q+1)$-fold suspension 
in Definition \ref{cofdef} we find that
$$S^{q(m+1)+1}\vee S^{(q+1)(m+1)} \to 
S^{(q+1)(m+1)}\vee S^{(q+1)(m+1)+m} \to
\Sigma^{q+1} C_q(S^m).$$
Any map $S^r\to S^s$ is null homotopic when $r<s$ and
a map $S^r \to S^r$ is determined, up to homotopy, 
by its degree. By Proposition \ref{Euler} we have that
the degree of the above self map of $S^{(q+1)(m+1)}$
is $\chi (S^m)$ which is $0$ for m odd and $2$ for $m$ 
even. So $\Sigma^{q+1} C_q(S^m)$ is homotopy equivalent to
the following wedge of spheres and Moore spaces:
$$
\begin{cases}
S^{(q+1)(m+1)}\vee S^{(q+1)(m+1)+m}\vee S^{q(m+1)+2}\vee
S^{(q+1)(m+1)+1} &, m \text{ odd } \\
M(\ZZ /2 , (q+1)(m+1))\vee 
S^{(q+1)(m+1)+m}\vee S^{q(m+1)+2} &, m \text{ even. }
\end{cases}
$$
Thus $\Sigma^{(m-2)(q+1)} C_q(S^m)$ is stably equivalent to
$$
\begin{cases}
S^{2(q+1)(m-1)}\vee S^{2(q+1)(m-1)+m}\vee S^{(2q+1)(m-1)}\vee
S^{2(q+1)(m-1)+1} &, m \text{ odd } \\
M(\ZZ /2 , 2(q+1)(m-1))\vee 
S^{2(q+1)(m-1)+m}\vee S^{(2q+1)(m-1)} &, m \text{ even. }
\end{cases}
$$
So we have determined the stable homotopy type of $\lm (S^m)$.

On the other hand the results in \cite{BM}, \cite{CC} give a
stable decomposition of $\Lambda \Sigma Y$ where $Y$ is a based
connected space. Let $C_k \subseteq S^1$ denote the cyclic group 
of order $k$ and put $D_k(Y)=S^1_+ \wedge_{C_k} Y^{\wedge k}$.
Then there is a stable equivalence 
$$\Lambda \Sigma Y \simeq \bigvee_{k\geq 1} D_k(Y).$$ 
We use this result for $Y=S^{m-1}$. Here each $D_k(S^{m-1})$ 
decomposes further.

The cofiber of the inclusion $i:C_{k+} \subseteq S^1_+$ is a 
$k$-fold wedge of circles. We have an associated Barrat-Puppe sequence
$$
\begin{CD}
C_{k+} \wedge_{C_k} Y^{\wedge k} @>i\wedge id>> D_k(Y) @>>> 
(S^1\wedge C_{k+})\wedge_{C_k} Y^{\wedge k} @>>> \\
@| @| @| \\
Y^{\wedge k} @>>> D_k(Y) @>>> S^1\wedge Y^{\wedge k} @>>>
\end{CD}
$$
So for $Y=S^{m-1}$ we find 
$$S^{(m-1)k} \to D_k(S^{m-1}) \to S^{(m-1)k+1} 
\xrightarrow{\partial_{m-1}} S^{(m-1)k+1} \to 
\Sigma D_k(S^{m-1}) \to \dots$$
By analyzing the Barrat-Puppe sequence

$$C_{k+} \xrightarrow{i} S^1_+ \to 
S^1\wedge C_{k+}\xrightarrow{\partial} S^1\wedge C_{k+}\to S^1\wedge
S^1_+ \to \dots $$
via the mapping cone of $i$ we see that $(\partial_{m-1})_*= id - T_*$ 
on reduced integral homology. Shifting to integral cohomology
we see that the degree of $\partial_{m-1}$ is $2$ when 
$m-1$ is odd and $k-1$ is odd and $0$ otherwise. Thus
$$
\Sigma D_k(S^{m-1}) \simeq 
\begin{cases}
M(\ZZ /2 , (m-1)k+1) &, \text{ for $m$ and $k$ even} \\
S^{(m-1)k+1}\vee S^{(m-1)k+2} &, \text{ otherwise.} 
\end{cases}
$$
So $D_k(S^{m-1})$ is stably equivalent to $M(\ZZ /2, (m-1)k)$
when $m$ and $k$ are both even and stably equivalent to 
$S^{(m-1)k}\vee S^{(m-1)k+1}$ otherwise. The result follows.
\end{proof}

\section{Appendix: A property of the shuffle map}

Let $A\simp$ and $B\simp$ be simplicial modules over a 
commutative ring $k$. Recall that the shuffle map 
$D: C_*(A)\otimes C_*(B) \to C_*(A\otimes B)$
is defined as follows \cite{ML} VIII.8:
$$D_n (a\otimes b) = \sum_{(\mu , \nu)} (-1)^{\epsilon (\mu )}
s_{\nu_q} \dots s_{\nu_1}a \otimes s_{\mu_p} \dots s_{\mu_1}b$$
where $a\in A_p$ and $b\in B_q$ with $p+q=n$. The sum is 
taken over all $(p,q)$ shuffles $(\mu , \nu )$ cf. \cite{ML}
and the sign is given by the formula
\begin{equation}
\label{eq:sign}
\epsilon (\mu )=\sum_{i=1}^p (\mu_i - (i-1)).
\end{equation}

We want to investigate how the shuffle map commutes with 
the simplicial boundary maps. Recall that
\begin{equation}
\label{eq:identity}
d_is_j= \begin{cases}
s_{j-1}d_i &, i<j \\
id &, i=j \text{ or } i=j+1 \\
s_jd_{i-1} &, i>j+1.
\end{cases}
\end{equation}

\begin{proposition}
\label{EZ}
There is a commutative diagram as follows:
$$
\begin{CD}
(C_*(A) \otimes C_*(B))_n @>D_n>> C_n(A\otimes B) \\
@V\partial_0 VV @Vd_0\otimes d_0 VV \\
(C_*(A)\otimes C_*(B))_{n-1} @>D_{n-1}>> C_{n-1}(A\otimes B)
\end{CD}
$$
where the map $\partial_0$ is defined by
$$\partial_0: A_p\otimes B_q \to (A_{p-1}\otimes B_q) \oplus 
(A_p \otimes B_{q-1}); \quad 
\partial_0=d_0\otimes id + (-1)^p id\otimes d_0.$$
Here $A_j=B_j=0$ when $j<0$ by convention.
When $i\geq 1$ there are $k$-linear maps  
$F_t, G_t:(C_*(A)\otimes C_*(B))_{n-1} \to C_{n-1}(A\otimes B)$ 
for $1\leq t \leq i$ such that 
$$
(d_i \otimes d_i)\circ D_n=
\sum_{t=1}^i 
\big( F_t\circ (d_t\otimes id)+G_t\circ (id\otimes d_t)\big) .
$$
\end{proposition}

\begin{proof}
We rewrite the sum defining $(d_0\otimes d_0)\circ D_n$
by commuting $d_0$ and $s_{\mu_j}$, using formula 
(\ref{eq:identity}). There are two possibilities:
$$
d_0 s_{\mu_p}\dots s_{\mu_1}=
\begin{cases}
s_{\mu_p-1}\dots s_{\mu_1-1}d_0 &\text{ if } \mu_1 \geq 1\\
s_{\mu_p-1}\dots s_{\mu_2-1} & \text{ if } \mu_1 =0
\end{cases}
$$
corresponding to the two cases
\[
\begin{split}
d_0 s_{\nu_q}\dots & s_{\nu_1}\otimes 
d_0 s_{\mu_p}\dots s_{\mu_1}=\\
& \begin{cases}
s_{\nu_q-1}\dots s_{\nu_2-1}\otimes 
s_{\mu_p-1}\dots s_{\mu_1-1} d_0
& \text{ if } \nu_1=0 \\
s_{\nu_q-1}\dots s_{\nu_1-1} d_0\otimes 
s_{\mu_p-1}\dots s_{\mu_2-1}
& \text{ if } \mu_1=0.
\end{cases}
\end{split}
\]

To the shuffle $(\mu,\nu)$ of the elements in $[n-1]$, 
we associate a new shuffle $(\mu-1,\nu-1)$ on $[n-2]$ by the rule
$(\mu-1)_i=\mu_i-1$ and $(\nu-1)_i=\nu_i-1$. There is of course 
the problem that $(\mu-1)_1$ can be negative. Actually, 
either $(\mu-1)_1=-1$ or $(\nu-1)_1=-1$, and we delete 
the corresponding element. What we have left is a genuine 
shuffle on $[n-2]$.

This defines a bijection between shuffles on $[n-1]$ 
such that $\mu_1=0$ and shuffles on $[n-2]$. 
Similarly, it defines a bijection between shuffles 
of $[n-1]$ element such that $\nu_1=0$ and shuffles on 
$[n-2]$ elements. 

The sign of the shuffle $(\mu-1,\nu-1)$ is, according to
formula (\ref{eq:sign}) given by
$$
\epsilon(\mu-1)=
\begin{cases}
\sum_{1\leq i\leq p-1} ((\mu-1)_i-(i-1))=\epsilon(\mu)
&\text{ if } \mu_1=0\\
\sum_{1\leq i\leq p} ((\mu-1)_i-(i-1))=\epsilon(\mu)-p
&\text{ if } \nu_1=0.  
\end{cases}
$$ 
So we can express the sum defining $(d_0\otimes d_0)\circ D_n$
as a sum of shuffles (with sign) applied to either
$(d_0\otimes id)$ or $(id\otimes d_0)$. 
Thus the diagram commutes.

Let $i\geq 1$. We use formula (\ref{eq:identity}) to rewrite 
$d_i s_{\mu_p}\dots s_{\mu_1}$, moving the $d_i$ toward the 
right. There are two possibilities. 

If $\{ i,i-1\} \cap \{ \mu_1, \dots ,\mu_p \} = \emptyset$ 
we get that the element equals some 
$s_{\mu_p^\prime}\dots s_{\mu_1^\prime}d_r$. 
Looking closer at the
formula (\ref{eq:identity}), we see that since $i\geq 1$, 
we have that $i\geq r \geq 1$. 
This means that there are maps $\alpha$ and $\beta$
such that 
$$(d_i\otimes d_i)(s_{\nu_q}\dots s_{\nu_1} a \otimes 
s_{\mu_p}\dots s_{\mu_1}b)= 
(\alpha \otimes \beta)\circ (id \otimes d_r)(a\otimes b).$$ 

Similarly, if 
$\{ i,i-1\} \cap \{ \nu_1 , \dots , \nu_q \} = \emptyset$ 
there are maps $\alpha^\prime$ and $\beta^\prime$ and
an integer $s$ with $i\geq s\geq 1$ such that
$$(d_i \otimes d_i) (s_{\nu_q}\dots s_{\nu_1}a\otimes 
s_{\mu_p}\dots s_{\mu_1}b)=
(\alpha^\prime \otimes \beta^\prime )\circ 
(d_s \otimes id)(a\otimes b).$$

We now only have to take care of the summands of $D_n$
corresponding to shuffles such that either $i=\mu_r$ 
and $i-1=\nu_s$ for some $r,s$ or $i-1 =\mu_r$ and 
$i=\nu_s$ for some $r,s$. 
These shuffles come in pairs:
If $(\mu,\nu)$ is a shuffle with $\mu_r=i$ and $\nu_s=i-1$, then
we can define a different shuffle $(\mu^\prime, \nu^\prime)$ by
interchanging $\mu_r$ and $\nu_s$. Thus
$$
\mu^\prime_t=
\begin{cases}
\mu_t & \text{ if } t\neq r\\
i-1 & \text{ if } t=r.
\end{cases}
$$
Then $(\mu,\nu)$ and $(\mu^\prime,\nu^\prime)$ will according to
formula (\ref{eq:sign}) contribute equal
terms with opposite signs in the sum computing 
$(d_i\otimes d_i)\circ D_n$.
These terms cancel, and the second part of the proposition is proved.
\end{proof}


\begin{thebibliography}{WWW}

\bibitem[A]{A} J. F. Adams, On the non-existence of elements of 
Hopf invariant one, Ann. of Math. 72 (1960), 20-104.

\bibitem[BK]{BK} A.K. Bousfield \& D.M. Kan, Homotopy limits and
localizations, Springer Lecture Notes in Mathematics, Vol 304(1972)

\bibitem[BM]{BM} C.-F. B\" odigheimer \& I. Madsen, 
Homotopy quotients of mapping spaces and their stable splitting,
Quart. J. Math. Oxford (2), 39 (1988), 401-409.


\bibitem[BO1]{BO2} M. B\"okstedt \& I. Ottosen, A spectral sequence
for string cohomology, Copenhagen preprint (2002).

\bibitem[BO2]{BO3} M. B\"okstedt \& I. Ottosen, An alternative 
approach to homotopy operations, \AA rhus preprint (2004).

\bibitem[CC]{CC} G. Carlsson \& R. Cohen, The cyclic groups and the
free loop space, Comment. Math. Helv. 62 (1987), 423-449.

\bibitem[CS]{CS} M. Chas \& D. Sullivan, String Topology, 
to appear in Annals of Math.

\bibitem[D]{D} W. G. Dwyer, Homotopy operations for simplicial 
commutative algebras, Trans A.M.S. 260 (1980), 421-435.

\bibitem[G]{G} P. G. Goerss, Andr\' e-Quillen cohomology and the 
homotopy groups of mapping spaces: Understanding the $E_2$-term of
the Bousfield-Kan spectral sequence, JPAA 63 (1990), 113-153.

\bibitem[KY]{KY} K. Kuribayashi \& T. Yamaguchi, The cohomology
algebra of certain free loop spaces, Fund. Math. 154 (1997),
57-73.

\bibitem[ML]{ML} S. MacLane, Homology, Springer 1975.

\bibitem[MS]{MS} J. W. Milnor \& J. D. Stasheff, Characteristic
Classes, Annals of Mathematics Studies, Princeton University Press
1974.

\bibitem[T]{T} H. Toda, Note on cohomology ring of certain spaces,
Proc. Amer. Math. Soc. 14 (1963), 89-95.

\bibitem[Z]{Z} W. Ziller, The free loop space of globally symmetric 
spaces, Invent. Math. 41 (1977), 1-22.

\end{thebibliography}
\end{document}